\documentclass[number,dvips]{arxbj}
\usepackage{stfloats}
\usepackage{graphicx}

\aid{0}
\volume{15}
\issue{1}
\pubyear{2009}
\firstpage{146}
\lastpage{177}
\doi{10.3150/08-BEJ142}

\makeatletter
\fnbelowfloat
\newtheorem{theorem}{Theorem}
\newtheorem{prop}{Proposition}
\newtheorem{lemme}{Lemma}
\newremark{Remarque}{Remark}
\makeatother

\begin{document}
\begin{frontmatter}

\title{Approximation of the distribution of a stationary Markov
process with application to option pricing}
\runtitle{Approximation of the distribution of a stationary Markov
process}

\begin{aug}
\author[a]{\fnms{Gilles} \snm{Pag\`{e}s}\thanksref{a}\ead[label=e1]{gpa@ccr.jussieu.fr}\corref{}} \and
\author[b]{\fnms{Fabien} \snm{Panloup}\thanksref{b}\ead[label=e2]{fpanloup@insa-toulouse.fr}}
\runauthor{G. Pag\`{e}s and F. Panloup}
\pdfauthor{Gilles Pages and Fabien Panloup}
\address[a]{Laboratoire de Probabilit\'es et Mod\`eles Al\'eatoires,
UMR 7599, Universit\'e Paris 6, Case 188, 4 pl. Jussieu, F-75252 Paris
Cedex 5. \printead{e1}}
\address[b]{Laboratoire de Statistiques et Probabilit\'{e}s,
Universit\'{e} Paul Sabatier \&
INSA Toulouse, 135, Avenue de Rangueil, 31077 Toulouse Cedex 4.
\printead{e2}}
\end{aug}

\received{\smonth{4} \syear{2007}}
\revised{\smonth{3} \syear{2008}}

%
\begin{abstract}
We build a sequence of empirical measures on the space
$\mathbb{D}(\mathbb{R}_+,\mathbb{R}^d)$ of $\mathbb{R}^d$-valued
cadlag functions on
$\mathbb{R}_+$ in order to approximate the law of a stationary
$\mathbb{R}^d$-valued Markov and Feller process $(X_t)$. We obtain some
general results on the convergence of this sequence. We then apply
them to Brownian diffusions and solutions to L\'{e}vy-driven
SDE's under some Lyapunov-type stability assumptions.
As a numerical application of this work, we show that this
procedure provides an efficient means of option pricing in stochastic
volatility models.
\end{abstract}

%
\pdfkeywords{Euler scheme, Levy process,
numerical approximation, option pricing, stationary process,
stochastic volatility model, tempered stable process,}

\end{frontmatter}

\section{Introduction}\label{s1}
\subsection{Objectives and motivations}

In this paper, we deal with an $\mathbb{R}^d$-valued Feller Markov
process $(X_t)$ with semigroup $(P_t)_{t\ge0}$ and assume
that $(X_t)$ admits an invariant distribution ${\nu_0}$. The aim
of this work is to propose a way to approximate the whole
stationary distribution $\mathbb{P}_{\nu_0}$ of $(X_t)$. More
precisely, we want to construct a sequence of weighted occupation
measures $(\nu^{(n)}(\omega,\mathrm{d}\alpha))_{n\ge1}$ on the
Skorokhod space
$\mathbb{D}(\mathbb{R}_+,\mathbb{R}^d)$ such that
$\nu^{(n)}(\omega,F)\stackrel{n\rightarrow+\infty}{\longrightarrow
}\int F(\alpha)
\mathbb{P}_{\nu_0}(\mathrm{d}\alpha)$ a.s.
for a class of functionals $F\dvtx\mathbb{D}(\mathbb{R}_+,\mathbb{R}^d)$
which includes bounded
continuous functionals for the Skorokhod topology.

One of our motivations is to develop a new numerical method for
option pricing in stationary stochastic volatility models which are
slight modifications of the classical stochastic volatility
models, where we suppose that the volatility evolves under its
stationary regime.

\subsection{Background and construction of the procedure}

This work follows on from a series of recent papers due to
Lamberton and Pag\`{e}s (\cite{LP1,LP2}), Lemaire
(\cite{lemaire2,lemaire1}) and Panloup
(\cite{panloup1,panloup3,panloup2}), where the problem of the
approximation of the invariant distribution is investigated for
Brownian diffusions and for L\'{e}vy-driven SDE's.
\footnote{Note that
computing the invariant distribution is equivalent to computing
the marginal laws of the stationary process $(X_t)$ since ${\nu_0}
P_t={\nu_0}$ for every $t\ge0$.}
In these papers, the algorithm
is based on an adapted Euler scheme with decreasing step
$(\gamma_k)_{k\ge1}$. To be precise, let $(\Gamma_n)$ be the
sequence of discretization times: $\Gamma_0=0$,
$\Gamma_n=\sum_{k=1}^n\gamma_k$ for every $n\ge1$, and assume that
$\Gamma_n\rightarrow+\infty$ when $n\rightarrow+\infty$. Let
$(\bar{X}_{\Gamma_n})_{n\ge0}$ be the Euler scheme obtained by
``freezing'' the coefficients between the $\Gamma_n$'s and let
$(\eta_n)_{n\ge1}$ be a sequence of positive weights such that
$H_n:=\sum_{k=1}^n\eta_k\rightarrow+\infty$ when
$k\rightarrow+\infty$. Then, under some Lyapunov-type stability
assumptions adapted to the stochastic processes of interest, one
shows that for a large class of steps and weights
$(\eta_n,\gamma_n)_{n\ge1}$,
%
%
\begin{equation}\label{marginalconvergence}
\bar{\nu}_n(\omega,f):=\frac{1}{H_n}\sum_{k=1}^n
\eta_kf(\bar{X}_{\Gamma_{k-1}})
\stackrel{n\rightarrow+\infty}{\longrightarrow}\int f(x){\nu
_0}(\mathrm{d}x)
\qquad\mbox{a.s.},
\end{equation}
(at least)
\footnote{The class of functions for which
(\ref{marginalconvergence})
holds depends on the stability of the dynamical system. In
particular, in the Brownian diffusion case, the convergence may
hold for continuous functions with subexponential growth, whereas
the class of functions strongly depends on the moments of the L\'{e}vy
process when the stochastic process is a L\'{e}vy-driven SDE.}
for every bounded continuous function $f$.

Since the problem of the approximation of the invariant
distribution has been deeply studied for a wide class of Markov
processes (Brownian diffusions and L\'{e}vy-driven SDE's) and since
the proof of (\ref{marginalconvergence}) can be adapted to other
classes of Markov processes under some specific Lyapunov
assumptions, we choose in this paper to consider a general Markov
process and to assume the existence of a time discretization scheme
$(\bar{X}_{\Gamma_k})_{k\ge0}$ such that
(\ref{marginalconvergence}) holds for the class of bounded
continuous functions. The aim of this paper is then to investigate
the convergence properties of a functional
version of the sequence
$(\bar{\nu}_n(\omega,\mathrm{d}\alpha))_{n\ge1}$.

Let $(X_t)$ be a Markov and Feller process and let
$(\bar{X}_t)_{t\ge0}$ be a stepwise constant time discretization
scheme of
$(X_t)$ with non-increasing step sequence $(\gamma_n)_{n\ge1}$
satisfying
%
%
\begin{equation}\label{nonincreagam}
\lim_{n\rightarrow+\infty}\gamma_n=0,\qquad
\Gamma_n:=\sum_{k=1}^n\gamma_k\stackrel{n\rightarrow+\infty
}{\longrightarrow}
+ \infty.
\end{equation}
Letting $\Gamma_0:=0$ and $\bar{X}_0=x_0\in\mathbb{R}^d$, we assume that
%
%
\begin{equation}\label{contschemerel}
\bar{X}_t=\bar{X}_{\Gamma_n}\qquad
\forall t\in[\Gamma_n,\Gamma_{n+1}[
\end{equation}
and that
$(\bar{X}_{\Gamma_n})_{n\ge0}$ can be simulated recursively.

We denote by $(\mathcal{F}_t)_{t\ge0}$ and
$(\bar{\mathcal{F}}_t)_{t\ge0}$ the usual augmentations of the natural
filtrations $(\sigma(X_s,0\le s\le t))_{t\ge0}$ and
$(\sigma(\bar{X}_s,0\le s\le t))_{t\ge0}$, respectively.

For $k\ge0$, we denote by $(\bar{X}_t^{(k)})_{t\ge0}$ the
shifted process defined by
\[
\bar{X}_t^{(k)}:=\bar{X}_{\Gamma_k+t}.
\]
In particular,
$\bar{X}_{t}^{(0)}=\bar{X}_{t}$. We define a sequence of random
probabilities $(\nu^{(n)}(\omega,\mathrm{d}\alpha))_{n\ge1}$ on
$\mathbb{D}(\mathbb{R}_+,\mathbb{R}^d)$ by
\[
\nu^{(n)}(\omega,\mathrm{d}\alpha)=\frac{1}{H_n}\sum_{k=1}^n\eta_k
{\bf{1}}_{{\{\bar{X}}^{({k-1})}(\omega)\in\mathrm{d}\alpha\}},
\]
where
$(\eta_k)_{k\ge1}$ is a sequence of weights. For $t\ge0$,
$(\nu_t^{(n)}(\omega,\mathrm{d}x))_{n\ge1}$ will denote the
sequence of
``marginal'' empirical measures on $\mathbb{R}^d$ defined by
\[
\nu^{(n)}_t(\omega,\mathrm{d}x)=\frac{1}{H_n}\sum_{k=1}^n\eta_k
{\bf{1}}_{\{\bar{X}^{({k-1})}_t(\omega)\in\mathrm{d}x\}}.
\]

\subsection{Simulation of $(\nu^{(n)}(\omega,F))_{n\ge1}$}\label
{simfuncpropos}

For every functional
$F\dvtx\mathbb{D}(\mathbb{R}_+,\mathbb{R}^d)\rightarrow\mathbb
{R}$, the
following recurrence relation holds for every \mbox{$n\ge1$}:
%
%
\begin{equation}\label{foncrecurrelation}
{\nu^{(n+1)}}(\omega,F)
= {\nu^{(n)}}(\omega,F)+\frac{\eta_{n+1}}{H_{n+1}}
\bigl(F\bigl(X^{(n)}(\omega)\bigr)-{\nu^{(n)}}(\omega,F)\bigr).
\end{equation}
Then, if $T$ is a positive number and
$F\dvtx\mathbb{D}(\mathbb{R}_+,\mathbb{R}^d)\rightarrow\mathbb{R}$
is a functional depending only on
the trajectory between 0 and $T$, $({\nu^{(n)}}(\omega,F))_{n\ge1}$ can
be simulated by the following procedure.

\textbf{Step 0.} (i) Simulate $(\bar{X}^{(0)}_t)_{t\ge0}$ on
$[0,T]$, that is, simulate $(\bar{X}_{\Gamma_k})_{k\ge0}$ for
$k=0,\ldots,{{N}(0,T)}$, where
%
%
\begin{eqnarray}\label{foncNnt}
{N}(n,T) :\!\!&=&\inf\{k\ge n, \Gamma_{k+1}-\Gamma_n> T\}
\nonumber\\[-8pt]
\\[-8pt]
&=& \max\{k\ge0,\Gamma_k-\Gamma_n\le T\},\qquad n\ge0, T>0.
\nonumber
\end{eqnarray}
Note that $n\mapsto N(n,t)$ is an increasing sequence since
$(\gamma_n)$ is non-increasing, and that
\[
\Gamma_{N(n,T)}-\Gamma_n\le T<\Gamma_{N(n,T)+1}-\Gamma_n.
\]

(ii) Compute $F((\bar{X}_t^{(0)})_{t\ge0})$ and
$\nu^{(1)}(\omega,F)$. Store the values
of $(\bar{X}_{\Gamma_k})$ for $k=1,\ldots,{N}(0,T)$.

\textbf{Step} $\bolds{n}$ \textbf{(}$\bolds{n\ge1}$\textbf{).} (i) Since the
values $(\bar{X}_{\Gamma_k})_{k\ge0}$ are stored for
$k=n,\ldots,N(n-1,T)$, simulate $(\bar{X}_{\Gamma_k})_{k\ge0}$
for $k={N}(n-1,T)+1,\ldots,{N}(n,T)$ in order to obtain a path of
$(\bar{X}_t^{(n)})$ on $[0,T]$.

(ii) Compute $F((\bar{X}_t^{(n)})_{t\ge0})$ and use
(\ref{foncrecurrelation}) to compute ${\nu^{(n+1)}}(\omega,F)$.
Store the values of $(\bar{X}_{\Gamma_k})$ for
$k=n+1,\ldots,{N}(n,T)$.
\begin{Remarque}
As shown in the description of the procedure, one
generally has to store the vector
$[\bar{X}_{\Gamma_n},\ldots,\bar{X}_{\Gamma_{N(n,T)}}]$ at time $n$.
Since $(\gamma_n)$ is a sequence with infinite sum that decreases to
0, it follows that the size of this vector
increases ``slowly'' to $+\infty$. For instance,
if $\gamma_n=C n^{-\rho}$ with $\rho\in(0,1)$, its size is of order
$n^{\rho}$. However,
it is important to remark that even though the number of values to
be stored tends to $+\infty$, that is not always the case
for the number of operations at each step.
Indeed, since $\bar{X}^{(n+1)}$ is obtained by shifting
$\bar{X}^{(n)}$, it is usually
possible to use, at step $n+1$, the preceding computations and to
simulate the sequence $(F(\bar{X}^{(n)}))_{n\ge0}$ in a
``quasi-recursive'' way. For instance, such remark holds for Asian options
because the associated pay-off can be expressed as a function of an
additive functional (see Section \ref{foncfinance} for simulations).
\end{Remarque}

Before outlining the sequel of the paper, we list some
notation linked to the spaces $\mathbb{D}(\mathbb{R}_+,\mathbb
{R}^d)$ and
$\mathbb{D}([0,T],\mathbb{R}^d)$ of cadlag
$\mathbb{R}^d$-valued functions on $\mathbb{R}_+$ and $[0,T]$,
respectively, endowed with the Skorokhod topology.
First, we denote by $d_1$ the Skorokhod distance on
$\mathbb{D}([0,1],\mathbb{R}^d)$ defined for
every $\alpha$, $\beta\in\mathbb{D}([0,1],\mathbb{R}^d)$ by
\[
d_1(\alpha,\beta)=\inf_{\lambda\in\Lambda_1}
\biggl\{\max\biggl(\sup_{t\in[0,1]}|\alpha(t)-\beta(\lambda(t))|,
\sup_{0\le s<t\le1}
\biggl|\log\biggl(\frac{\lambda(t)-\lambda(s)}{t-s}\biggr)
\biggr| \biggr) \biggr\},
\]
where $\Lambda_1$ denotes the set of increasing homeomorphisms of
$[0,1]$. Second, for $T>0$,
$\phi_T\dvtx\mathbb{D}(\mathbb{R}_+,\mathbb{R}^d)\mapsto\mathbb
{D}([0,1],\mathbb{R}^d)$
is the function defined by $(\phi_T(\alpha))(s)=\alpha({sT})$
for every $s\in[0,1]$. We then denote by $d$ the distance on
$\mathbb{D}(\mathbb{R}_+,\mathbb{R}^d)$ defined for every
$\alpha, \beta\in\mathbb{D}(\mathbb{R}_+,\mathbb{R}^d)$ by
%
%
\begin{equation}\label{distancecomplete}
d(\alpha,\beta)=\int_0^{+\infty} \mathrm{e}^{-t}
\bigl(1\wedge d_1(\phi_t(\alpha),\phi_t(\beta)) \bigr)\,\mathrm{d}t.
\end{equation}
We recall that $(\mathbb{D}(\mathbb{R}_+,\mathbb{R}^d),d)$ is a
Polish space and that the induced
topology is the usual Skorokhod topology on
$\mathbb{D}(\mathbb{R}_+,\mathbb{R}^d)$ (see, e.g., Pag\`es
\cite{pagesthese}). For every $T>0$, we set
\[
\mathcal{D}_T=\bigcap_{s>T}\sigma(\pi_u,0\le u\le s),
\]
where $\pi_s\dvtx
\mathbb{D}(\mathbb{R}_+,\mathbb{R}^d)\rightarrow\mathbb{R}^d$
is defined by $\pi_s(\alpha)=\alpha(s)$. For a functional
$F\dvtx\mathbb{D}(\mathbb{R}_+,\mathbb{R}^d)\rightarrow\mathbb{R}$,
$F_T$ denotes the functional defined for every
$\alpha\in\mathbb{D}(\mathbb{R}_+,\mathbb{R}^d)$ by
%
%
\begin{equation}\label{functionalt}
F_T(\alpha)=F(\alpha^T)\qquad\mbox{with }
\alpha^T(t)=\alpha(t\wedge T)\ \forall t\ge0.
\end{equation}
Finally, we will say that a
functional $F\dvtx\mathbb{D}(\mathbb{R}_+,\mathbb{R}^d)\rightarrow
\mathbb{R}$ is
$Sk$-continuous if $F$ is continuous for the Skorokhod topology on
$\mathbb{D}(\mathbb{R}_+,\mathbb{R}^d)$ and the notation
``$\stackrel{(Sk)}{\Longrightarrow}$'' will denote the weak
convergence on $\mathbb{D}(\mathbb{R}_+,\mathbb{R}^d)$.

In Section \ref{generalresultsfunc}, we state our main results for a
general $\mathbb{R}^d$-valued
Feller Markov process.
Then, in Section \ref{appliSDE}, we apply them to Brownian diffusions
and L\'{e}vy-driven SDE's.
Section \ref{proofsfunc} is devoted to the proofs of the main general
results. Finally, in Section
\ref{foncfinance}, we complete this
paper with an application to option pricing in stationary stochastic
volatility models.

\section{General results}\label{generalresultsfunc}
In this section, we state the results on convergence of the
sequence $(\nu^{(n)}(\omega,\mathrm{d}\alpha))_{n\ge1}$ when
$(X_t)$ is a
general Feller Markov process.

\subsection{Weak convergence to the stationary regime}

As explained in the \hyperref[s1]{Introduction}, since the a.s. convergence of
$(\nu^{(n)}_0(\omega,\mathrm{d}x))_{n\ge1}$
to the invariant distribution $\nu_0$ has already been deeply
studied for a large class of Markov processes (Brownian diffusions
and L\'{e}vy driven SDE's), our approach will be to derive the
convergence of $(\nu^{(n)}(\omega,\mathrm{d}\alpha))_{n\ge1}$ toward
$\mathbb{P}_{\nu_0}$ from that of $(\nu^{(n)}_0(\omega,\mathrm
{d}x))_{n\ge1}$ to
the invariant distribution ${\nu_0}$. More precisely, we will
assume in Theorem \ref{prinfonc1} that

$\mathbf{(C_{0,1})}$: $(X_t)$ admits a unique invariant
distribution ${\nu_0}$ and
\[
\nu^{(n)}_0(\omega,\mathrm{d}x)\stackrel{n\rightarrow+\infty
}{\Longrightarrow}\nu_0(\mathrm{d}x)\qquad\mbox{a.s.},
\]
whereas in Theorem \ref{prinfonc2}, we will only assume that

$\mathbf{(C_{0,2})}$: $(\nu_0^{(n)}(\omega,\mathrm{d}x))_{n\ge1}$
is a.s. tight on $\mathbb{R}^d$.

We also introduce three other assumptions, $\mathbf{(C_1)}$,
$\mathbf{(C_2)}$ and $\mathbf{(C_{3,\varepsilon})}$, regarding the
continuity in probability of the flow $x\mapsto(X_t^x)$, the asymptotic
convergence of the shifted time discretization scheme
to the true process $(X_t)$ and the steps and weights, respectively.

$\mathbf{(C_1)}$: For every $x_0\in\mathbb{R}^d$,
$\epsilon>0$
and $T>0$,
%
%
\begin{equation}\label{djioe2}
\limsup_{x_0\rightarrow x}\mathbb{P}\biggl(\sup_{0\le t\le
T}|X_t^x-X_t^{x_0}|\ge\epsilon\biggr)=0.
\end{equation}

$\mathbf{(C_2)}$: $(\bar{X}_t)$ is a non-homogeneous
Markov process and for every $n\ge0$, it is possible to construct
a family of stochastic processes
$({{Y}}_t^{(n,x)})_{x\in\mathbb{R}^d}$ such that
\begin{longlist}
\item$\mathcal{L}({{Y}}^{(n,x)})\stackrel{\mathbb{D}(\mathbb{R}_+,
\mathbb{R}^d)}{=}\mathcal{L}(\bar{X}^{(n)}|\bar{X}^{(n)}_0=x)$;
\item for every compact set $K$ of $\mathbb{R}^d$, for every $T\ge0$,
%
\begin{equation}\label{djioe3}
\sup_{x\in K}\sup_{0\le t\le T}
\bigl|{{Y}}^{(n,x)}_t-X_t^x\bigr|
\stackrel{n\rightarrow+\infty}{\longrightarrow}0\qquad\mbox{in
probability}.
\end{equation}
\end{longlist}
$\mathbf{(C_{3,\varepsilon})}$: For every $n\ge1$,
$\eta_n\le C\gamma_n H_n^{\varepsilon}$.
\begin{Remarque}
Assumption $\mathbf{(C_2)}$ implies, in particular, that
asymptotically and uniformly on compact sets of $\mathbb{R}^d$, the law
of the approximate process $(\bar{X}^{(n)})$, given
its initial value, is close to that of the true process.

If there exists a unique invariant distribution ${\nu_0}$, the
second part of $\mathbf{(C_2)}$ can be relaxed to the following, less
stringent, assertion: for all $\epsilon>0$, there exists a compact set
$A_\epsilon\subset\mathbb{R}^d$ such that
${\nu_0}(A_\epsilon^c)\le\epsilon$ and such that
%
%
\begin{equation}\label{djioe}
\sup_{x\in A_\epsilon}\sup_{0\le t\le T}
\bigl|{{Y}}^{(n,x)}_t-X_t^x\bigr|
\stackrel{n\rightarrow+\infty}{\longrightarrow}0\qquad\mbox{in
probability}.
\end{equation}
This weaker assumption can some times be needed
in stochastic volatility models like the Heston model (see
Section \ref{foncfinance} for details).
\end{Remarque}

The preceding assumptions are all that we require for the
convergence of $(\nu^{(n)}(\omega,\mathrm{d}\alpha))_{n\ge1}$ to
$\mathbb{P}_{{\nu}_0}$ along the bounded $Sk$-continuous functionals,
that is, for the a.s. weak convergence on
$\mathbb{D}(\mathbb{R}_+,\mathbb{R}^d)$. However, the integration of
non-bounded continuous functionals
$F\dvtx{\mathbb D}([0,T],\mathbb{R}^d)\rightarrow\mathbb{R}$
will need some additional assumptions, depending on the stability
of the time discretization scheme and on the steps and weights
sequences. We
will suppose that $F$ is dominated (in a sense to be specified
later) by a function $\mathcal{V}:\mathbb{R}^d\rightarrow\mathbb{R}_+$
that satisfies the following assumptions for some $s\ge2$ and
$\varepsilon<1$.

$\mathbf{H(s,\varepsilon)}$: For every $T>0$,
\begin{longlist}
\item$\displaystyle\sup_{n\ge1}\mathbb{E}
\biggl[\sup_{0\le t\le T}\mathcal{V}^s\bigl({{Y}}^{(n,x)}_t\bigr
)\biggr]
\le C_T\mathcal{V}^s(x)$,\vspace*{2pt}
\item$\displaystyle\sup_{n\ge1}\nu_0^{(n)}(\mathcal{V})<+\infty
$,\vspace*{2pt}
\item$\displaystyle\sum_{k\ge1}\frac{\eta_k}{H_k^2}\mathbb
{E}[\mathcal{V}^2
(\bar{X}_{\Gamma_{k-1}})]<+\infty$,\vspace*{2pt}
\item$\displaystyle\sum_{k\ge1}\frac{\Delta
N(k,T)}{H_k^{s}}\mathbb{E}
\bigl[\mathcal{V}^{s(1-\varepsilon)}(\bar{X}_{\Gamma_{k-1}})\bigr
]<+\infty$,
\end{longlist}
where $T\mapsto C_T$ is locally bounded on $\mathbb{R}_+$
and $\Delta N(k,T)=N(k,T)-N(k-1,T)$.

For every $\varepsilon<1$, we then set
\[
\mathcal{K}(\varepsilon)= \{\mathcal{V}\in\mathcal{C}(\mathbb{R}^d,
\mathbb{R}_+),\mathbf{H(s,\varepsilon)} \mbox{ holds for some }s\ge
2 \}.
\]
\begin{Remarque}
Apart from assumption (i), which is a classical condition
on the finite time horizon control, the
assumptions in $\mathbf{H(s,\varepsilon)}$ strongly rely on the stability
of the time discretization scheme (and then, to that of the true
process). More precisely, we will
see when we apply our general results to SDE's that these
properties are some consequences of the Lyapunov assumptions
needed for the tightness of $(\nu_0^{(n)}(\omega,\mathrm{d}x))_{n\ge1}$.
\end{Remarque}

We can now state our first main result.
\begin{theorem}\label{prinfonc1}
Assume $\mathbf{(C_{0,1})}$, $\mathbf{(C_1)}$, $\mathbf{(C_2)}$
and $\mathbf{(C_{3,\varepsilon})}$ with
$\varepsilon\in(-\infty,1)$. Then, a.s., for every bounded
$Sk$-continuous functional
$F\dvtx\mathbb{D}(\mathbb{R}_+,\mathbb{R}^d)\rightarrow\mathbb{R}$,
%
%
\begin{equation}\label{conlip}
{\nu}^{(n)}(\omega,F)\stackrel{n\rightarrow+\infty
}{\longrightarrow}\int
F(\alpha)\mathbb{P}_{\nu_0}(\mathrm{d}\alpha),
\end{equation}
where $\mathbb{P}_{\nu_0}$ denotes the stationary distribution of
$(X_t)$ (with initial law ${\nu_0}$).

Furthermore, for every $T>0$, for every non-bounded $Sk$-continuous
functional
$F\dvtx\mathbb{D}(\mathbb{R}_+,\mathbb{R}^d)\rightarrow
\mathbb{R}$,
(\ref{conlip}) holds a.s. for $F_T$ (defined by (\ref{functionalt})) if
there exists $\mathcal{V}\in\mathcal{K}(\varepsilon)$ and $\rho\in[0,1)$
such that
%
%
\begin{equation}\label{convnonbounded}
|F_T(\alpha)|\le C\sup_{0\le t\le T}\mathcal{V}^\rho(\alpha_t)
\qquad\forall \alpha\in\mathbb{D}(\mathbb{R}_+,\mathbb{R}^d).
\end{equation}
\end{theorem}

In the second result, the uniqueness of the invariant
distribution is not required and the sequence
$(\nu^{(n)}_0(\omega,\mathrm{d}x))_{n\ge1}$ is only supposed to be
tight.
\begin{theorem}\label{prinfonc2}
Assume $\mathbf{(C_{0,2})}$, $\mathbf{(C_1)}$, $\mathbf{(C_2)}$
and $\mathbf{(C_{3,\varepsilon})}$ with
$\varepsilon\in(-\infty,1)$. Assume that
$(\nu^{(n)}_0(\omega,\mathrm{d}x))_{n\ge1}$ is a.s. tight on
$\mathbb{R}^d$. We then have the following.
\begin{longlist}
\item The sequence $(\nu^{(n)}(\omega,\mathrm{d}\alpha))_{n\ge1}$
is a.s. tight
on $\mathbb{D}(\mathbb{R}_+,\mathbb{R}^d)$ and a.s., for every convergent
subsequence $(n_k(\omega))_{n\ge1}$, for every bounded
$Sk$-continuous functional
$F\dvtx\mathbb{D}(\mathbb{R}_+,\mathbb{R}^d)\rightarrow\mathbb{R}$,
%
%
\begin{equation}\label{conv2functheo}
\nu^{(n_k(\omega))}(\omega,F)
\stackrel{n\rightarrow+\infty}{\longrightarrow}
\int F(\alpha)\mathbb{P}_{{\nu}_\infty}(\mathrm{d}\alpha),
\end{equation}
where $\mathbb{P}_{\nu_\infty}$ is the law of $(X_t)$ with initial law
$\nu_\infty$ being a weak limits for
$(\nu^{(n)}_0(\omega,\mathrm{d}x))_{n\ge1}$.

Furthermore, for
every $T>0$, for every non-bounded $Sk$-continuous functional
$F\dvtx\mathbb{D}(\mathbb{R}_+,\mathbb{R}^d)\rightarrow
\mathbb{R}$,
(\ref{conv2functheo}) holds a.s. for $F_T$ if
(\ref{convnonbounded}) is satisfied with $\mathcal{V}\in\mathcal
{K}(\varepsilon)$ and
$\rho\in[0,1)$.
\item If, moreover,
%
%
\begin{equation}\label{invarianceautomatique}
\frac{1}{H_n}\sum_{k=1}^n\max_{l\ge k+1}\frac{|\Delta
\eta_\ell|}{\gamma_{\ell}}\stackrel{n\to +\infty}{\longrightarrow}0,
\end{equation}
then ${\nu}_\infty$ is necessarily an invariant
distribution for the Markov process $(X_t)$.
\end{longlist}
\end{theorem}
\begin{Remarque}
Condition (\ref{invarianceautomatique}) holds for a
large class of steps and weights.
For instance, if $\eta_n=C_1n^{-\rho_1}$ and
$\gamma_n=C_2n^{-\rho_2}$ with $\rho_1\in[0,1]$ and
$\rho_2\in(0,1]$, then (\ref{invarianceautomatique}) is satisfied
if $\rho_1=0$ or if $\rho_1\in(\max(0,2\rho_2-1),1)$.
\end{Remarque}

\subsection{Extension to the non-stationary case}\label{remarqueabs}

Even though the main interest of this algorithm is the
weak approximation of the process when stationary, we observe that
when $\nu_0$ is known, the algorithm can be used to approximate
$\mathbb{P}_{\mu_0}$ if $\mu_0$ is a probability on $\mathbb{R}^d$ that
is absolutely continuous with respect to $\nu_0$.

Indeed, assume that $\mu_0(\mathrm{d}x)=\phi(x)\nu_0(\mathrm
{d}x)$, where
$\phi\dvtx\mathbb{R}^d\rightarrow\mathbb{R}$ is a continuous
non-negative function. For a functional
$F\dvtx\mathbb{D}(\mathbb{R}_+,\mathbb{R}^d)\rightarrow\mathbb{R}$,
denote by $F^\phi$ the functional defined on
$\mathbb{D}(\mathbb{R}_+,\mathbb{R}^d)$ by
$F^{\phi}(\alpha)=F(\alpha)\phi(\alpha(0))$.

Then, if $\nu^{(n)}(\omega,\mathrm{d}\alpha)\stackrel
{(Sk)}{\Rightarrow}
\mathbb{P}_{\nu_0}(\mathrm{d}\alpha)$ a.s., we also have the
following convergence: a.s., for every bounded $Sk$-continuous
functional $F\dvtx\mathbb{D}(\mathbb{R}_+,\mathbb{R}^d)\rightarrow
\mathbb{R}$,
\[
{\nu}^{(n)}(\omega,F^\phi)\stackrel{n\rightarrow+\infty
}{\longrightarrow}
\int F^\phi(\alpha)\mathbb{P}_{\nu_0}(\mathrm{d}\alpha)
=\int F(\alpha)\mathbb{P}_{\mu_0}(\mathrm{d}\alpha).
\]

\section{Application to Brownian diffusions and L\'{e}vy-driven
SDE's}\label{appliSDE}

Let $(X_t)_{t\ge0}$ be a cadlag stochastic process solution to the SDE
%
%
\begin{equation}\label{edss}
\mathrm{d}X_{t}=b(X_{t^{-}})\,\mathrm{d}t
+ \sigma(X_{t^{-}})\,\mathrm{d}W_{t}+\kappa(X_{t^{-}})
\,\mathrm{d}Z_{t},
\end{equation}
where $b\dvtx\mathbb{R}^d\rightarrow\mathbb{R}^d$,
$\sigma\dvtx\mathbb{R}^d\mapsto\mathbb{M}_{d,\ell}$
(set of $d\times\ell$ real matrices) and
$\kappa\dvtx\mathbb{R}^d\mapsto\mathbb{M}_{d,\ell}$ are
continuous functions
with sublinear growth, $(W_{t})_{t\ge0}$ is an $\ell$-dimensional
Brownian motion and $(Z_{t})_{t\ge0}$ is an integrable purely discontinuous
$\mathbb{R}^{\ell}$-valued L\'{e}vy process independent of
$(W_{t})_{t\ge0}$
with L\'{e}vy measure $\pi$ and characteristic function given for
every $t\ge0$ by
\[
\mathbb{E}[\mathrm{e}^{\mathrm{i}\langle u,Z_t\rangle}]=\exp\biggl[t \biggl(\int
\mathrm{e}^{\mathrm{i}\langle u,y\rangle}-1
- \mathrm{i}\langle u,y\rangle\pi(\mathrm{d}y) \biggr)
\biggr].
\]
Let $(\gamma_n)_{n\ge1}$ be a non-increasing step sequence
satisfying (\ref{nonincreagam}). Let $(U_n)_{n\ge1}$ be a sequence
of i.i.d. random variables such that
$U_1\stackrel{\mathcal{L}}{=}\mathcal{N}(0,I_{\ell})$ and let
$\xi:=(\xi_n)_{n\ge1}$ be a sequence of
independent $\mathbb{R}^{\ell}$-valued random variables, independent of
$(U_n)_{n\ge1}$. We then denote by $(\bar{X}_t)_{t\ge0}$ the
stepwise constant Euler scheme of $(X_t)$ for which
$(\bar{X}_{{\Gamma_{n}}})_{n\ge0}$ is recursively defined by
$\bar{X}_0=x\in\mathbb{R}^d$ and
%
%
\begin{equation}\label{Eulerscheme}
{\bar{X}}_{{\Gamma_{n+1}}}=\bar{X}^{}_{{\Gamma_{n}}}
+\gamma_{n+1} b(\bar{X}^{}_{{\Gamma_{n}}})
+ \sqrt{\gamma_{n+1}}\sigma(\bar{X}^{}_{{\Gamma_{n}}})U_{n+1}
+\kappa(\bar{X}^{}_{{\Gamma_{n}}})\xi_{n+1}.
\end{equation}
We recall that the increments of $(Z_t)$ cannot be simulated in
general. That is why we generally need to construct the sequence
$(\xi_n)$ with some approximations of the true increments. We will
come back to this construction in Section \ref{levydriven}.

As in the general case, we denote by $(\bar{X}^{(k)})_{k\ge0}$ and
$(\nu^{(n)}(\omega,\mathrm{d}\alpha))_{n\ge1}$ the sequences of associated
shifted Euler schemes and empirical measures, respectively.

Let us now introduce some Lyapunov assumptions for the
SDE. Let $\mathcal{E}\mathcal{Q}(\mathbb{R}^d)$ denote the set of
\textit{essentially quadratic} $\mathcal{C}^2$-functions
$V\dvtx\mathbb{R}^d\rightarrow\mathbb{R}_+^*$ such that
$\lim V(x)=+\infty$ as
$|x|\rightarrow+\infty$, $|\nabla V|\le C \sqrt{V}$
and $D^2V$
is bounded. Let $a\in(0,1]$ denote the mean reversion intensity.
The Lyapunov (or mean reversion) assumption is the following.

$\mathbf{(S_a)}$:
There exists a function $V\in\mathcal{E}\mathcal{Q}(\mathbb{R}^d)$
such that:
\begin{longlist}
\item\,\ $|b|^2\le CV^a$,
$\operatorname{Tr}(\sigma\sigma^*(x))+\|\kappa(x)\|^2
\stackrel{|x|\rightarrow+\infty}{=}o(V^a(x))$;
\item\,\ there exist $\beta\in\mathbb{R}$ and $\rho>0$
such that $\langle\nabla V,b\rangle\le\beta-\rho V^a$.
\end{longlist}
From now on, we separate the Brownian diffusions and L\'{e}vy-driven
SDE cases.

\subsection{Application to Brownian diffusions}

In this part, we assume that $\kappa=0$. We recall a result by
Lamberton and
Pag\`es \cite{LP2}.
\begin{prop}\label{funcpremiere1}
Let $a\in(0,1]$ such that
$\mathbf{(S_{a})}$ holds. Assume that the sequence
$(\eta_n/\gamma_n)_{n\ge1}$ is non-increasing.
\begin{longlist}[(b)]
\item[(a)] Let $(\theta_n)_{n\ge1}$ be a sequence of positive
numbers such that $\sum_{n\ge1}\theta_n\gamma_n<+\infty$ and that
there exists $n_0\in\mathbb{N}$ such that $(\theta_n)_{n\ge n_0}$ is
non-increasing. 
Then, for every positive $r$,
\[
\sum_{n\ge1}\theta_n\gamma_n\mathbb{E}
[V^{r}(\bar{X}_{\Gamma_{n-1}})] <+\infty.
\]
\item[(b)] For every $r>0$,
%
%
\begin{equation}\label{Vsuptendu}
\sup_{n\ge1}{\nu}_0^{(n)}(\omega,V^r)<+\infty\qquad\mbox{a.s.}
\end{equation}
Hence, the sequence $(\nu^{(n)}_0(\omega,\mathrm{d}x))_{n\ge1}$ is a.s.
tight.
\item[(c)] Moreover, every weak limit of this sequence is an
invariant probability for the SDE (\ref{edss}). In particular, if
$(X_t)_{t\ge0}$ admits a unique invariant probability ${\nu_0}$, then
for every continuous function~$f$ such that $f\le C V^r$ with
$r>0$, $\lim_{n\rightarrow\infty} \nu_0^{(n)}(\omega,f)={\nu
_0}(f)$ a.s.
\end{longlist}
\end{prop}
\begin{Remarque} 
For instance, if $V(x)=1+|x|^2$, then the preceding convergence holds
for every continuous function with polynomial growth.
According to Theorem 3.2 in Lemaire \cite{lemaire2}, it is possible to
extend these results to continuous functions with exponential
growth, but it then strongly depends on $\sigma$. Further the
conditions on steps and weights can be less restrictive and may
contain the case $\eta_n=1$, for instance (see Remark 4 of Lamberton
and Pag\`es \cite{LP2} and Lemaire \cite{lemaire2}).
\end{Remarque}

We then derive the following result from the preceding
proposition and from Theorems \ref{prinfonc1}~and~\ref{prinfonc2}.
\begin{theorem}\label{prinfonc4} 
Assume that $b$ and $\sigma$ are locally
Lipschitz functions and that $\kappa=0$. Let $a\in(0,1]$ such
that $\mathbf{(S_{a})}$ holds and assume that $(\eta_n/\gamma_n)$
is non-increasing.
\begin{longlist}[(b)]
\item[(a)] The sequence $(\nu^{(n)}(\omega,\mathrm{d}\alpha
))_{n\ge1}$ is
a.s. tight on $\mathcal{C}(\mathbb{R}_+,\mathbb{R}^d)$
\footnote{$\mathcal{C}(\mathbb{R}_+,\mathbb{R}^d)$
denotes the space of continuous
functions on $\mathbb{R}_+$ with values in
$\mathbb{R}^d$ endowed with the topology of uniform convergence on
compact sets.}
and every weak limit of
$(\nu^{(n)}(\omega,\mathrm{d}\alpha))_{n\ge1}$ is the distribution
of a
stationary process solution to (\ref{edss}). In
particular, when uniqueness holds for the invariant distribution
$\nu_0$, a.s., for every bounded continuous functional
$F\dvtx\mathcal{C}(\mathbb{R}_+,\mathbb{R}^d)\rightarrow\mathbb{R}$,
%
%
\begin{equation}\label{converbrowniansde}
\nu^{(n)}(\omega,F)\stackrel{n\rightarrow+\infty}{\longrightarrow}
\int F(x)\mathbb{P}_{\nu_0}(\mathrm{d}x).
\end{equation}
\item[(b)] Furthermore, if there exists $s\in(2,+\infty)$ and
$n_0\in \mathbb{N}$ such that
%
%
\begin{equation}\label{condipaspoids}
\biggl(\frac{\Delta N(k,T)}{\gamma_k H_k^s} \biggr)_{n\ge n_0}
\mbox{ is non-increasing and }
\sum_{k\ge1}\frac{\Delta N(k,T)}{H_k^s}<+\infty,
\end{equation}
then, for every $T>0$, for every non-bounded continuous functional
$F\dvtx\mathcal{C}(\mathbb{R}_+,\mathbb{R}^d)\rightarrow\mathbb{R}$,
(\ref{converbrowniansde}) holds
for $F_T$ if the following condition is satisfied:
\[
\exists r>0 \quad\mbox{such that}\quad
|F_T(\alpha)|\le C\sup_{0\le t\le T}V^r(\alpha_t)
\qquad\forall\alpha\in\mathcal{C}(\mathbb{R}_+,\mathbb{R}^d).
\]
\end{longlist}
\end{theorem}
\begin{Remarque}\label{criterepas}
If $\eta_n=C_1n^{-\rho_1}$ and
$\gamma_n=C_2n^{-\rho_2}$ with $0<\rho_2\le\rho_1\le1$,
then for $s\in(1,+\infty)$, (\ref{condipaspoids}) is fulfilled if
and only if $s>1/(1-\rho_1)$. It follows that there exists $s\in
(2,+\infty)$ such that
(\ref{condipaspoids}) holds as soon as $\rho_1<1$.
\end{Remarque}
\begin{pf*}{Proof of Theorem \ref{prinfonc4}}
We want to apply Theorem \ref{prinfonc2}. First, by Proposition
\ref{funcpremiere1}, assumption $\mathbf{(C_{0,2})}$ is fulfilled
and every weak limit of $(\nu_0^{(n)}(\omega,\mathrm{d}x))$ is an invariant
distribution. Second, it is well known that $\mathbf{(C_1)}$ and
$\mathbf{(C_2)}$ are fulfilled when $b$ and $\sigma$ are locally
Lispchitz sublinear functions. Then, since
$\mathbf{(C_{3,\varepsilon})}$ holds with $\varepsilon=0$,
(\ref{converbrowniansde}) holds for every bounded continuous
functional $F$. Finally, one checks that $\mathbf{H(s,0)}$ holds
with $\mathcal{V}:=V^r$ ($r>0$). It is classical that assumption (a)
is true when $b$ and $\sigma$ are sublinear. Assumption (b)
follows from Proposition \ref{funcpremiere1}(b). Let
$\theta_{n,1}=\eta_n/(\gamma_nH_n^2)$ and
$\theta_{n,2}=\Delta N(n,T)/(\gamma_n H_n^s)$. Using (\ref{condipaspoids})
and the fact that $(\eta_n/\gamma_n)$ is non-increasing yields that
$(\theta_{n,1})$ and $(\theta_{n,2})$ satisfy the conditions of
Proposition \ref{funcpremiere1} (see (\ref{preuvedec}) for
details). Then, (iii) and (iv) of $\mathbf{H(s,0)}$ are
consequences of Proposition \ref{funcpremiere1}(a). This
completes the proof.
\end{pf*}

\subsection{Application to L\'{e}vy-driven SDE's}\label{levydriven}

When we want to extend the results obtained for Brownian SDE's to
L\'{e}vy-driven SDE's, one of the main difficulties comes from the
moments of the jump component (see Panloup \cite{panloup1} for details).
For simplification, we assume here that $(Z_t)$ has a moment of
order $2p\ge2$, that is, that its L\'{e}vy measure $\pi$ satisfies the
following assumption with $p\ge1$:
\[
{\bf{(H^1_p)}}\dvtx \int_{|y|>1}\pi(\mathrm{d}y)|y|^{2p}<+\infty.
\]
We also introduce an assumption about the behavior of the moments
of the L\'{e}vy measure at 0:
\[
\qquad{\bf{(H^2_q)}}\dvtx\int_{|y|\le1}\pi(\mathrm{d}y)|y|^{2q}
< +\infty,\qquad q\in[0,1].
\]
This assumption ensures that $(Z_t)$ has finite $2q$-variations.
Since $\int_{|y|\le1}|y|^2\pi(\mathrm{d}y)$ is finite, this is
always satisfied for $q=1$.

Let us now specify the law of $(\xi_n)$ introduced in
(\ref{Eulerscheme}). When the increments of $(Z_t)$ can be exactly
simulated, we denote by (E) the Euler scheme and by $(\xi_{n,E})$
the associated sequence
\[
\xi_{n,E}\stackrel{\mathcal{L}}{=}Z_{\gamma_n}\qquad\forall n\ge1.
\]
\\
When the increments of $(Z_t)$ cannot be simulated, we introduce some
\textit{approximated} Euler schemes (P) and (W)
built with some sequences $(\xi_{n,P})$ and $(\xi_{n,W})$ of
approximations of the true increment (see Panloup \cite{panloup3} for
more detailed presentations of these schemes).

In scheme (P),
\[
\xi_{n,P}\stackrel{\mathcal{L}}{=}Z_{\gamma_n,n},
\]
where $(Z_{\bolds{\cdot},n})_{n\ge1}$
a sequence of compensated compound Poisson processes obtained by
truncating the small jumps of $(Z_t)_{t\ge0}$:
%
%
\begin{equation}\label{sautapprox}
Z_{t,n}:=\sum_{0<s\le t}\Delta Z_s 1_{\{|\Delta Z_s|>u_n \}}
- t\int_{|y|>u_n}y\pi(\mathrm{d}y)\qquad\forall t\ge0,
\end{equation}
where $(u_n)_{n\ge1}$ is a sequence of positive numbers such that
$u_n\rightarrow0$. We recall that $Z_{\bolds{\cdot},n}
\stackrel{n\rightarrow+\infty}{\longrightarrow} Z$ locally
uniformly in $L^2$ (see, e.g., Protter \cite{protter}).

As shown in Panloup \cite{panloup3}, the error induced by this
approximation is very large when the local behavior of the small
jumps component is irregular. However, it is possible to refine
this approximation by a \textit{Wienerization} of the small jumps, that
is, by replacing the small jumps by a
linear transform of a Brownian motion instead of discarding them (see
Asmussen and Rosinski~\cite{rosinski}). The
corresponding scheme is denoted by (W) with $\xi_{n,W}$ satisfying
\[
\xi_{n,W}\stackrel{\mathcal{L}}{=}\xi_{n,P}
+ \sqrt{\gamma_n}Q_n\Lambda_n\qquad\forall n\ge1,
\]
where $(\Lambda_n)_{n\ge
1}$ is a sequence of i.i.d. random variables, independent of
$(\xi_{n,P})_{n\ge1}$ and $(U_n)_{n\ge1}$, such that
$\Lambda_1\stackrel{\mathcal{L}}{=}\mathcal{N}(0,I_\ell)$ and
$(Q_n)$ is a
sequence of $\ell\times\ell$ matrices such that
\[
(Q_nQ_n^*)_{i,j}=\int_{|y|\le u_k} y_i y_j\pi(\mathrm{d}y).
\]
We recall the following
result obtained in Panloup \cite{panloup1} in our slightly
simplified framework.
\begin{prop}\label{funcpremiere}
Let $a\in(0,1]$, $p\ge1$ and $q\in[0,1]$
such that $\mathbf{(H^1_p)}$, $\mathbf{(H^2_q)}$ and
$\mathbf{(S_{a})}$ hold. Assume that the sequence
$(\eta_n/\gamma_n)_{n\ge1}$ is non-increasing. Then, the
following assertions hold for schemes \textup{(E)}, \textup{(P)} and
\textup{(W)}.
\begin{longlist}[(b)]
\item[(a)] Let $(\theta_n)$ satisfy the conditions of Proposition
\ref{funcpremiere1}. Then,
$\sum_{n\ge1}\theta_n\gamma_n\mathbb{E}[V^{p+a-1}(\bar
{X}_{\Gamma_{n-1}})]<+\infty$.
\item[(b)] We have
%
%
\begin{equation}\label{Vsuptendu2}
\sup_{n\ge1}{\nu}_0^{(n)}
(\omega,V^{{p}/{2}+a-1})<+\infty\qquad\mbox{a.s.}
\end{equation}
Hence, the sequence $(\nu^{(n)}_0(\omega,\mathrm{d}x))_{n\ge1}$ is a.s.
tight as soon as $p/2+a-1>0$.
\item[(c)] Moreover, if ${\mathrm{Tr}}(\sigma\sigma^*)+\|\kappa\|
^{2q}\le C
V^{{p}/{2}+a-1}$, then every weak limit of this sequence is
an invariant probability for the SDE (\ref{edss}). In particular,
if $(X_t)_{t\ge0}$ admits a unique invariant probability
${\nu_0}$, for every continuous function $f$ such that
$f=o(V^{{p}/{2}+a-1})$,
$\lim_{n\rightarrow\infty}\nu_0^{(n)}(\omega,f)={\nu_0}(f)$ a.s.
\end{longlist}
\end{prop}
\begin{Remarque}
For schemes (E) and (P), the above proposition is a
direct consequence of Theorem~2
and Proposition 2 of Panloup \cite{panloup1}. As concerns scheme (W), a
straightforward adaptation of the proof yields the result.
\end{Remarque}

Our main functional result for L\'{e}vy-driven SDE's is
then the following.
\begin{theorem}\label{prinfonc3}
Let $a\in(0,1]$ and $p\ge1$ such that
$p/2+a-1>0$ and let $q\in[0,1]$. Assume
$\mathbf{(H_p^1)}$, $\mathbf{(H_q^2)}$ and $\mathbf{(S_{a})}$.
Assume that $b$, $\sigma$ and $\kappa$ are locally Lipschitz
functions. If, moreover, $(\eta_n/\gamma_n)_{n\ge1}$ is
non-increasing, then the following result holds for schemes
\textup{(E)}, \textup{(P)} and \textup{(W)}.
\begin{longlist}
\item[(a)] The sequence $(\nu^{(n)}(\omega,\mathrm{d}\alpha))_{n\ge1}$
is a.s. tight on $\mathbb{D}(\mathbb{R}_+,\mathbb{R}^d)$. Moreover, if
%
\begin{equation}\label{dnof}
\operatorname{Tr}(\sigma\sigma^*)+\|\kappa\|^{2q}
\le CV^{{p}/{2}+a-1}
\quad\mbox{or}\quad
\frac{1}{H_n} \sum_{k=1}^n\max_{l\ge k+1}
\frac{|\Delta\eta_\ell|}{\gamma_{\ell-1}}
\stackrel{n\to+\infty}{\longrightarrow} 0,
\end{equation}
then every weak limit of
$(\nu^{(n)}(\omega,\mathrm{d}\alpha))_{n\ge1}$ is the
distribution of a stationary process solution to~(\ref{edss}).
\item[(b)] Assume that the invariant distribution is unique.
Let $\varepsilon\le0$ such that $\mathbf{(C_{3,\varepsilon})}$
holds. Then, a.s.,
for every $T>0$, for every $Sk$-continuous functional
$F\dvtx\mathbb{D}(\mathbb{R}_+,\mathbb{R}^d)\rightarrow\mathbb{R}$,
(\ref{converbrowniansde}) holds for $F_T$ if there exist
$\rho\in[0,1)$ and $s\ge2$, such that
\[
|F_T(\alpha)|\le C\sup_{0\le t\le T}V^{(\rho(p+a-1))/{s}}
(\alpha_t)\qquad\forall\alpha\in\mathbb{D}(\mathbb{R}_+,\mathbb{R}^d)
\]
and if
%
%
\begin{equation}\label{condipoids2}
\biggl(\frac{\Delta N(k,T)}{\gamma_k
H_k^{s(1-\varepsilon)}} \biggr)_{n\ge n_0}
\mbox{ is non-increasing and }
\sum_{k\ge1}\frac{\Delta N(k,T)}{H_k^{s(1-\varepsilon)}}<+\infty.
\end{equation}
\end{longlist}
\end{theorem}
\begin{Remarque}
In (\ref{dnof}), both assumptions imply the invariance
of every weak limit of $(\nu_0^{(n)}(\omega,\mathrm{d}x))$. These two
assumptions are very different. The first is needed in
Proposition \ref{funcpremiere} for using the Echeverria--Weiss
invariance criteria (see Ethier and Kurtz \cite{bib4}, page~238,
Lamberton and Pag\`es \cite{LP1}
and Lemaire \cite{lemaire2}), whereas the second appears in Theorem
\ref{prinfonc2}, where our functional approach shows that under
some mild additional conditions on steps and
weights, every weak limit is always invariant.

For (\ref{condipoids2}), we refer to Remark \ref{criterepas} for
simple sufficient conditions when
$(\gamma_n)$ and $(\eta_n)$ are some polynomial steps and weights.
\end{Remarque}

\section[Proofs of Theorems 1 and 2]{Proofs of Theorems \protect\ref{prinfonc1}
and \protect\ref{prinfonc2}}\label{proofsfunc}

We begin the proof with some technical lemmas. In Lemma
\ref{caracterisation}, we show that the $a.s$ weak convergence of
the random measures $(\nu^{(n)}(\omega,\mathrm{d}\alpha))_{n\ge1}$
can be characterized by the convergence (\ref{conlip}) along the set of
bounded Lipschitz functionals $F$ for the distance $d$. Then, in
Lemma \ref{fundamental}, we show with some martingale arguments
that if the functional $F$ depends only on the restriction of the
trajectory to $[0,T]$, then the convergence of $(\nu^{(n)}(\omega,F))_{n\ge1}$ is equivalent
to that of a more regular sequence. This step is fundamental for
the sequel of the proof.

Finally, Lemma \ref{funclemme3} is needed for the proof of Theorem
\ref{prinfonc2}. We show that under some mild conditions on the
step and weight sequences, any Markovian weak limit of the
sequence $(\nu^{(n)}(\omega,\mathrm{d}\alpha))_{n\ge1}$ is stationary.

\subsection{Preliminary lemmas}

\begin{lemme}\label{caracterisation}
Let $(E,d)$ be a Polish space
and let $\mathcal{P}(E)$ denote the set of probability measures on
the Borel $\sigma$-field $\mathcal{B}(E)$, endowed with the weak
convergence topology.
Let $(\mu^{(n)}(\omega,\mathrm{d}\alpha))_{n\ge
1}$ be a sequence of random probabilities defined on $\Omega\times
\mathcal{B}(E)$.
\begin{longlist}[(b)]
\item[(a)] Assume that there exists $\mu^{(\infty)}\in\mathcal{P}(E)$ such
that for every bounded Lipschitz function $F\dvtx E\rightarrow\mathbb{R}$,
%
\begin{equation}\label{funcconv1}
\mu^{(n)}(\omega,F)\stackrel{n\rightarrow+\infty}
{\longrightarrow}\mu^{(\infty)}(F)\qquad \mbox{a.s.}
\end{equation}
Then, a.s., $(\mu^{(n)}(\omega,\mathrm{d}\alpha))_{n\ge1}$
converges weakly
to $\mu^{(\infty)}$ on $\mathcal{P}(E)$.
\item[(b)] Let $\mathcal{U}$ be a subset of $\mathcal{P}(E)$. Assume that for
every sequence $(F_{k})_{k\ge1}$ of Lipschitz and bounded
functions, a.s., for every subsequence
$(\mu^{(\phi_\omega(n))}(\omega,\mathrm{d}\alpha))$, there exists
a subsequence $(\mu^{(\phi_\omega\circ\psi_\omega(n))}(\omega,\mathrm{d}\alpha))$ and
a $\mathcal{U}$-valued random probability $\mu^{(\infty)}(\omega,\mathrm{d}\alpha)$
such that for every $k\ge1$,
%
\begin{equation}\label{funcconv2}
\mu^{(\psi_\omega\circ\phi_\omega(n))}(\omega,F_{k})
\stackrel{n\rightarrow+\infty}{\longrightarrow}
\mu^{(\infty)}(\omega,F_{k})\qquad\mbox{a.s.}
\end{equation}
Then, $(\mu^{(n)}(\omega,\mathrm{d}\alpha))_{n\ge1}$ is
a.s. tight with weak limits in $\mathcal{U}$.
\end{longlist}
\end{lemme}
\begin{pf}
We do not give a detailed proof of the next lemma, which is
essentially based
on the fact that in a separable metric space $(E,d)$, one can build a
sequence of bounded Lipschitz functions $(g_k)_{k\ge1}$
such that for any sequence $(\mu_n)_{n\ge1}$ of probability
measures on $\mathcal{B}(E)$, $(\mu_n)_{n\ge1}$ weakly converges
to a probability $\mu$ if and only if the convergence holds along the
functions $g_k,$ $k\ge1$ (see Parthasarathy \cite{partha}, Theorem
6.6, page~47 for a very similar result).
\end{pf}

For every $n\ge0$, for every $T>0$, we introduce
$\tau(n,T)$ defined by
%
\begin{equation}\label{fonctaunt}
\tau(n,T):=\min\{k\ge0, N(k,T)\ge n\}
= \min\{k\le n,\Gamma_k+T\ge\Gamma_n\}.
\end{equation}
Note that for $k\in\{0,\ldots,\tau(n,T)-1\}$,
$\{\bar{X}_t^{(k)},0\le t\le T\}$ is
$\mathcal{\bar{F}}_{\Gamma_n}$-measurable and
\[
T-\gamma_{\tau(n,T)-1}\le\Gamma_n-\Gamma_{\tau(n,T)}\le T.
\]
\begin{lemme} \label{fundamental}
Assume $\mathbf{(C_{3,\varepsilon})}$ with $\varepsilon<1$. Let
$F\dvtx \mathbb{D}(\mathbb{R}_+,\mathbb{R}^d)\rightarrow\mathbb{R}$ be
a $Sk$-continuous functional. Let $(\mathcal{G}_k)$ be a filtration
such that $\bar{\mathcal{F}}_{\Gamma_k}\subset\mathcal{G}_k$
for every $k\ge1$. Then, for any $T>0$:
\begin{longlist}[(b)]
\item[(a)] if $F_T$ (defined by (\ref{functionalt})) is bounded,
%
\begin{equation}\label{foncidfond}
\frac{1}{H_n}\sum_{k=1}^n \eta_k
\bigl(F_T\bigl(\bar{X}^{({k-1})}\bigr)-\mathbb{E}
\bigl[F_T\bigl(\bar{X}^{({k-1})}\bigr)/\mathcal{G}_{k-1}\bigr] \bigr)
\stackrel{n\rightarrow+\infty}{\longrightarrow}0\qquad \mbox{a.s.};
\end{equation}
\item[(b)] if $F_T$ is not bounded, (\ref{foncidfond}) holds if there exists
$\mathcal{V}\dvtx\mathbb{R}^d\rightarrow\mathbb{R}_+$, satisfying
$\mathbf{H(s,\varepsilon)}$ for some $s\ge2$, such that $|F_T(\alpha
)|\le C \sup_{0\le t\le T}\mathcal{V}(\alpha_t)$ for
every $\alpha\in\mathbb{D}(\mathbb{R}_+,\mathbb{R}^d)$; furthermore,
%
\begin{equation}\label{equiintfuncassump}
\sup_{n\ge1}\nu^{(n)}(\omega,F_T)<+\infty\qquad \mbox{a.s.}
\end{equation}
\end{longlist}
\end{lemme}
\begin{pf}
We prove (a) and (b) simultaneously. Let $\Upsilon^{(k)}$ be defined by
$\Upsilon^{(k)}=F_T(\bar{X}^{({k})})$. We have
%
\begin{eqnarray}
&& \frac{1}{H_n}\sum_{k=1}^n \eta_k
\bigl(\Upsilon^{(k-1)}-\mathbb{E}\bigl[\Upsilon^{(k-1)}/
\mathcal{G}_{k-1}\bigr] \bigr)
\nonumber\\
&&\quad = \frac{1}{H_n}\sum_{k=1}^n\eta_k
\bigl(\Upsilon^{(k-1)}-\mathbb{E}
\bigl[\Upsilon^{(k-1)}/\mathcal{G}_{n}\bigr] \bigr)
\label{foncdri27}\\
&&\qquad {} +\frac{1}{H_n}\sum_{k=1}^n \eta_k
\bigl(\mathbb{E}\bigl[\Upsilon^{(k-1)}/\mathcal{G}_{n}\bigr]
- \mathbb{E}\bigl[\Upsilon^{(k-1)}/\mathcal{G}_{k-1}\bigr] \bigr).
\label{foncdri28}
\end{eqnarray}
We have to prove that the right-hand side of (\ref{foncdri27}) and
(\ref{foncdri28}) tend to 0 a.s. when $n\rightarrow+\infty$.

We first focus on the right-hand side of (\ref{foncdri27}).
From the very definition of ${\tau(n,T)}$, we have that $\{\bar
{X}_t^{(k)},0\le t\le T\}$ is
$\bar{\mathcal{F}}_{\Gamma_n}$-measurable for
$k\in\{0,\ldots,\tau(n,T)-1\}$. Hence, since
$F_T$ is $\sigma(\pi_s,0\le s\le T)$-measurable and
$\bar{\mathcal{F}}_{\Gamma_n}\subset\mathcal{G}_{n}$, it follows that
$\Upsilon^{(k)}$ is $\mathcal{G}_{n}$-measurable and that
$\Upsilon^{(k)}=\mathbb{E}[\Upsilon^{(k)}/\mathcal{G}_{n}]$
for every $k\le\tau(n,T)-1$. Then, if $F_T$ is bounded,
we derive from $\mathbf{(C_{3,\varepsilon})}$ that
\begin{eqnarray*}
\Biggl|\frac{1}{H_n}\sum_{k=1}^n\eta_k
\bigl(\Upsilon^{(k-1)}-\mathbb{E}\bigl[\Upsilon^{(k-1)}\big/
\mathcal{G}_{n}\bigr] \bigr) \Biggr|
&\le& \frac{2\|F_T\|_{\mathrm{sup}}}{H_n}\sum_{k=\tau(n,T)+1}^n
\eta_k\le\frac{C}{H_n}\sum_{k=\tau(n,T)+1}^{n} \gamma_k H_k^{\varepsilon}
\\
&\le& \frac{C}{H_n^{1-\varepsilon}}
\bigl(\Gamma_n-\Gamma_{\tau(n,T)} \bigr)
\\
\\
&\le&\frac{C(T)}{H_n^{1-\varepsilon}}
\stackrel{n\rightarrow+\infty}{\longrightarrow}0\qquad\mbox{a.s.},
\end{eqnarray*}
where we used the fact that $(H_n)_{n\ge1}$ and $(\gamma_n)_{n\ge1}$ are
non-decreasing and non-increasing sequences, respectively.

Assume, now, that the assumptions of (b) are fulfilled with
$\mathcal{V}$ satisfying $\mathbf{H(s,\varepsilon)}$ for some $s\ge2$ and
$\varepsilon<1$. By the
Borel--Cantelli-like argument, it suffices to show that
%
\begin{equation}\label{foncdjks}
\sum_{n\ge1}\mathbb{E}
\Biggl[ \Biggl|\frac{1}{H_n^s}\sum_{k=\tau(n,T)+1}^n
\eta_k \bigl(\Upsilon^{(k-1)}-\mathbb{E}
\bigl[\Upsilon^{(k-1)}/\mathcal{G}_{n}\bigr] \bigr)
\Biggr|^s \Biggr]<+\infty.
\end{equation}
Let us prove (\ref{foncdjks}). Let $a_k:=\eta_k^{(s-1)/s}$ and
$b_k(\omega):=\eta_k^{{1}/{s}} (\Upsilon^{(k-1)}-\mathbb{E}
[\Upsilon^{(k-1)}/\mathcal{G}_{n}] )$. The H\"{o}lder inequality
applied with
$\bar{p}=s/(s-1)$ and $\bar{q}=s$ yields
\[
\Biggl|\sum_{k=\tau(n,T)+1}^n a_k b_k(\omega)\Biggr|^s
\le \Biggl(\sum_{k=\tau(n,T)+1}^n\eta_k \Biggr)^{s-1}
\Biggl(\sum_{k=\tau(n,T)+1}^n\eta_k
\bigl|\Upsilon^{(k-1)}-\mathbb{E}
\bigl[\Upsilon^{(k-1)}/\mathcal{G}_{n}\bigr] \bigr|^s \Biggr).
\]
Now, since $F_T(\alpha)\le
\sup_{0\le t\le T}\mathcal{V}(\alpha)$, it follows from the Markov
property and from $\mathbf{H(s,\varepsilon)}$(i) that
\[
\mathbb{E}\bigl[\bigl|F_T\bigl(\bar{X}^{(k)}\bigr)\bigr|^s
/\bar{\mathcal{F}}_{\Gamma_{k}}\bigr]
\le C \mathbb{E}\biggl[\sup_{0\le t\le T}
\mathcal{V}^s\bigl(\bar{X}_t^{(k)}\bigr)/\bar{\mathcal{F}}_{\Gamma_{k}}\biggr]
\le C_T \mathcal{V}^s(\bar{X}_{\Gamma_k}).
\]
Then, using the two
preceding inequalities and $\mathbf{(C_{3,\varepsilon})}$ yields
\begin{eqnarray*}
&& \mathbb{E} \Biggl[ \Biggl|\sum_{k=\tau(n,T)+1}^n\eta_k
\bigl(\Upsilon^{(k-1)}-\mathbb{E}
\bigl[\Upsilon^{(k-1)}/\mathcal{G}_{n}\bigr] \bigr) \Biggr|^s \Biggr]
\\
&&\quad \le  C \Biggl(\sum_{k=\tau(n,T)+1}^n\eta_k \Biggr)^{s-1}
\Biggl(\sum_{k=\tau(n,T)+1}^n\eta_k\mathbb{E}[\mathcal{V}^s(\bar{X}_{\Gamma_{k-1}})] \Biggr)
\\
&&\quad \le C \Biggl(\sum_{k=\tau(n,T)+1}^n\eta_k \Biggr)^{s}
\mathbb{E}\biggl[\sup_{k=\tau(n,T)+1}^n \mathcal{V}^s(\bar{X}_{\Gamma_{k-1}})\biggr]
\\
&&\quad \le C \Biggl(\sum_{k=\tau(n,T)+1}^n\gamma_kH_k^\varepsilon\Biggr)^{s}
\mathbb{E} \biggl[\sup_{t\in[0,S(n,T)]}
\mathcal{V}^s\bigl(\bar{X}^{\tau(n,T)}_{t}\bigr) \biggr],
\end{eqnarray*}
where $S(n,T)=\Gamma_{n-1}-\Gamma_{\tau(n,T)}$ and $C$ does not
depend $n$. By the definition of $\tau(n,T)$,
$S(n,T)\le T$. 
Then, again using $\mathbf{H(s,\varepsilon)}$(i) yields
\[
\sum_{n\ge1}\mathbb{E} \Biggl[\frac{1}{H_n^s}
\Biggl|\sum_{k=\tau(n,T)}^n\eta_k
\bigl(\Upsilon^{(k-1)}-\mathbb{E}\bigl[\Upsilon^{(k-1)}/
\mathcal{G}_{n}\bigr] \bigr) \Biggr|^s \Biggr]
\le C\sum_{n\ge 1}\frac{1}{H_n^{s(1-\varepsilon)}}
\mathbb{E}\bigl[\mathcal{V}^s\bigl(\bar{X}_{(\tau(n,T))}\bigr)\bigr].
\]
Since $n\mapsto N(n,T)$ is an increasing function,
$n\mapsto\tau(n,T)$ is a non-decreasing function and ${\rm
Card}\{n,\tau(n,T)=k\}=\Delta N(k+1,T):=N(k+1,T)-N(k,T)$. Then,
since $n\mapsto H_n$ increases, a change of variable yields
\begin{eqnarray*}
&& \sum_{n\ge1}\mathbb{E} \Biggl[\frac{1}{H_n^s} \Biggl|\sum_{k=\tau(n,T)+1}^n
\eta_k \bigl(\Upsilon^{(k-1)}-\mathbb{E}\bigl[\Upsilon^{(k-1)}/
\mathcal{G}_{n}\bigr] \bigr) \Biggr|^s \Biggr]
\\
&&\quad
\le C\sum_{k\ge1}\frac{\Delta N(k,T)}{H_k^{s(1-\varepsilon)}}
\mathbb{E}[\mathcal{V}^s(\bar{X}_{\Gamma_{k-1}})]<+\infty,
\end{eqnarray*}
by $\mathbf{H(s,\varepsilon)}$(iv).

Second, we prove that (\ref{foncdri28}) tends to 0.
For every $n\ge1$, we let
%
\begin{equation}
M_n=\sum_{k= 1}^n\frac{\eta_k}{H_k} \bigl(\mathbb{E}
\bigl[\Upsilon^{(k-1)}/\mathcal{G}_{n}\bigr]
- \mathbb{E}\bigl[\Upsilon^{(k-1)}/\mathcal{G}_{k-1}\bigr] \bigr).
\end{equation}
The process $(M_n)_{n\ge1}$ is a $(\mathcal{G}_{n})$-martingale and
we want to prove that this process is $L^2$-bounded. Set
$\Phi^{(k,n)}=\mathbb{E}[F_T(\bar{X}^{(k)})/\mathcal
{G}_{n}]-\mathbb{E}[F_T(\bar{X}^{(k)})/\mathcal{G}_k]$. Since
$F_T$ is $\sigma(\pi_s,0\le s\le T)$-measurable, the
random variable $\Phi^{(k,n)}$ is
$\bar{\mathcal{F}}_{\Gamma_{{N}(k,T)}}$-measurable. Then, for every
$i\in\{N(k,T),\ldots,n\}$,
$\Phi^{(k,n)}$ is \mbox{$\mathcal{G}_i$-measurable} so that
\[
\mathbb{E} \bigl[\Phi^{(i,n)}\Phi^{(k,n)} \bigr]
= \mathbb{E}\bigl[\Phi^{(k,n)}\mathbb{E}
\bigl[\Phi^{(i,n)}/\mathcal{G}_i\bigr]\bigr]=0.
\]
It follows that
%
\begin{equation}\label{rope5'}
\mathbb{E}[M_n^2]=\sum_{k\ge 1}\frac{\eta_k^2}{H_k^2}
\mathbb{E} \bigl[\bigl(\Phi^{({k-1},n)}\bigr)^2\bigr]
+ 2\sum_{k\ge 1}\frac{\eta_k}{H_k}\sum_{i=k+1}^{{N}(k-1,T)
\wedge n}\frac{\eta_i}{H_i}\mathbb{E}
\bigl[\Phi^{(i-1,n)}\Phi^{(k-1,n)} \bigr].
\end{equation}
Then,
%
\begin{eqnarray}\label{fonctionnelle168}
\sup_{n\ge1}\mathbb{E}[M_n^2] &\le& \sum_{k\ge 1}
\frac{\eta_k^2}{H_k^2}\sup_{n\geq1}
\mathbb{E} \bigl[\bigl(\Phi^{({k-1},n)}\bigr)^2 \bigr]
+ 2\sum_{k\ge1}\frac{\eta_k}{H_k}\sum_{i=k+1}^{{N}(k-1,T)}
\frac{\eta_i}{H_i}
\sup_{n\geq1}\mathbb{E}\bigl[\Phi^{(i-1,n)}\Phi^{(k-1,n)} \bigr]
\nonumber\\
&\le& C\Biggl(\sum_{k\ge 1}\frac{\eta_k}{H_k^{2-\varepsilon}}
\sup_{n\geq1}\mathbb{E}\bigl [\bigl(\Phi^{({k-1},n)}\bigr)^2 \bigr]
\\
&&\hspace*{14pt}{}
+ \sum_{k\ge1}\frac{\eta_k}{H_k^{2-\varepsilon}}
\sum_{i=k+1}^{{N}(k-1,T)}\gamma_i\sup_{n\geq1}\mathbb{E}
\bigl[\Phi^{(i-1,n)}\Phi^{(k-1,n)} \bigr]\Biggr),\qquad
\nonumber
\end{eqnarray}
where, in the second inequality, we used assumption
$\mathbf{(C_{3,\varepsilon})}$ and the decrease of
$i\mapsto1/H_i^{1-\varepsilon}$. Hence, if $F_T$ is bounded, using the
fact that $\sum_{i=k+1}^{{N}(k-1,T)}\gamma_i\le T$ yields
%
\begin{equation}\label{preuvedec}
\sup_{n\ge1}\mathbb{E}[M_n^2]
\le C\sum_{k\ge 1} \frac{\eta_k}{H_k^{2-\varepsilon}}
\le C \biggl(\frac{\eta_1}{H_1^{2-\varepsilon}}+\int_{\eta_1}^\infty
\frac{\mathrm{d}u}{u^{2-\varepsilon}} \biggr)<+\infty
\end{equation}
since $\varepsilon<1$. Assume, now, that the assumptions of (b) hold and
let $F_T$ be dominated by a function~$\mathcal{V}$ satisfying
$\mathbf{H(s,\varepsilon)}$. By the Markov property, the Jensen
inequality and $\mathbf{H(s,\varepsilon)}$(i),
\[
\mathbb{E} \bigl[\bigl(\Phi^{({k,n})}\bigr)^2\bigr]
\le C\mathbb{E}\biggl[\mathbb{E}\biggl[\sup_{0\le t\le T}
\mathcal{V}^2\bigl(\bar{X}_t^{(k)}\bigr)/
\bar{\mathcal{F}}_{\Gamma_k}\biggr]\biggr]
\le C_T \mathbb{E}[\mathcal{V}^2(\bar{X}_{\Gamma_k})].
\]
We then derive from the Cauchy--Schwarz inequality that for every
$n,k\ge1$, for every $i\in\{k,\ldots,N(k,T)\}$,
\begin{eqnarray*}
\bigl|\mathbb{E} \bigl[\Phi^{({i,n})}\Phi^{({k,n})}\bigr] \bigr|
\le C\sqrt{\mathbb{E}[\mathcal{V}^2(\bar{X}_{\Gamma_{i}})]}
\sqrt{\mathbb{E}[\mathcal{V}^2(\bar{X}_{\Gamma_{k}})]}
\le C\sup_{t\in[0,T]}\mathbb{E}
\bigl[\mathcal{V}^2\bigl(\bar{X}^{(k)}_{t}\bigr)\bigr]
\le C\mathbb{E}[\mathcal{V}^2(\bar{X}_{\Gamma_k})],
\end{eqnarray*}
where, in the last inequality, we once again used
$\mathbf{H(s,\varepsilon)}$(i). It follows that
\[
\sup_{n\ge1}\mathbb{E}[M_n^2]\le C\sum_{k\ge 1}
\frac{\eta_k}{H_k^{2-\varepsilon}}\mathbb{E}
[\mathcal{V}^2(\bar{X}_{\Gamma_{k-1}})]<+\infty,
\]
by $\mathbf{H(s,\varepsilon)}$(iii). Therefore,
(\ref{fonctionnelle168}) is finite and $(M_n)$ is bounded in
$L^2$. Finally, we derive from the Kronecker lemma that
\[
\frac{1}{H_n}\sum_{k= 1}^n\eta_k
\bigl(\mathbb{E}\bigl[F_T\bigl(\bar{X}^{({k-1})}\bigr)/
\mathcal{G}_{n}\bigr]-\mathbb{E}\bigl[F_T\bigl(\bar{X}^{({k-1})}\bigr)/
\mathcal{G}_{k-1}\bigr] \bigr)
\stackrel{n\rightarrow+\infty}{\longrightarrow}0\qquad \mbox{a.s.}
\]
As a consequence, $\sup_{n\ge1}\nu^{(n)}(\omega,F_T)<+\infty$ a.s. if
and only if
\[
\sup_{n\ge1}\frac{1}{H_n}\sum_{k=1}^n\mathbb{E}\bigl[F_T\bigl(\bar
{X}^{({k-1})}\bigr)/\mathcal{F}_{k-1}\bigr]<+\infty\qquad \mbox{a.s.}
\]
This last property is easily
derived from $\mathbf{H(s,\varepsilon)}$(i) and (ii). This
completes the proof.
\end{pf}
\begin{lemme} \label{lemmeflow}
\textup{(a)} Assume $\mathbf{(C_1)}$ and let $x_0\in\mathbb{R}^d$. We
then have $\lim_{x\rightarrow x_0}\mathbb{E}[d(X^x,X^{x_0})]=0$. In
particular, for every bounded Lispchitz (w.r.t. the distance $d$) functional
$F\dvtx\mathbb{D}(\mathbb{R}_+,\mathbb{R}^d)\rightarrow\mathbb{R}$,
the function $\Phi^F$ defined by
$\Phi^F(x)=\mathbb{E}[F(X^x)]$ is a (bounded) continuous function on
$\mathbb{R}^d$.

\textup{(b)} Assume $\mathbf{(C_2)}$. For every compact set $K\subset
\mathbb{R}^d$,
%
\begin{equation}\label{ucpe2}
\sup_{x\in K}\mathbb{E}[d({{Y}}^{n,x},X^{x})]
\stackrel{n\rightarrow+\infty}{\longrightarrow}0.
\end{equation}
Set $\Phi^F_{n}(x)=\mathbb{E}[F({{Y}}^{n,x})]$. Then, for every
bounded Lispchitz {functional }\mbox{
$F\dvtx\mathbb{D}(\mathbb{R}_+,\mathbb{R}^d)\rightarrow\mathbb{R}$,}
%
%
\begin{equation}\label{foncucp3}
\sup_{x\in K}|\Phi^F(x)-\Phi^F_{n}(x)|
\stackrel{n\rightarrow+\infty}{\longrightarrow}0
\qquad\mbox{for every compact set $K\subset\mathbb{R}^d$.}
\end{equation}
\end{lemme}
\begin{pf}
(a) By the definition of $d$, for every $\alpha$,
$\beta\in\mathbb{D}(\mathbb{R}_+,\mathbb{R}^d)$ and for every $T>0$,
%
%
\begin{equation}\label{ucpe}
d(\alpha,\beta)\le\biggl(1\wedge\sup_{0\le t\le T}
|\alpha(t)-\beta(t)| \biggr)+\mathrm{e}^{-T}.
\end{equation}
It easily
follows from assumption $\mathbf{(C_1)}$ and from the dominated
convergence theorem that
\[
\limsup_{x\rightarrow x_0}\mathbb{E}[d(X^x,X^{x_0})]
\le \mathrm{e}^{-T}\qquad\mbox{for every $T>0$.}
\]
Letting
$T\rightarrow+\infty$ implies that
$\lim_{x\rightarrow x_0}\mathbb{E}[d(X^x,X^{x_0})]=0$.

(b) We deduce from (\ref{ucpe}) and from assumption
$\mathbf{(C_2)}$ that for every compact set
$K \subset\mathbb{R}^d$, for every $T>0$,
\[
\limsup_{n\rightarrow+\infty}
\sup_{x\in K}\mathbb{E}
[d({{Y}}^{n,x},X^x)]\le\mathrm{e}^{-T}.
\]
Letting $T\rightarrow+\infty$ yields (\ref{ucpe2}).
\end{pf}
\begin{lemme}\label{funclemme3}
Assume that $(\eta_n)_{n\ge1}$ and $(\gamma_n)$ satisfy
$\mathbf{(C_{3,\varepsilon})}$ with $\varepsilon<1$ and
(\ref{invarianceautomatique}). Then:

\textup{(i)} for every $t\ge0$, for every bounded continuous function
$f\dvtx\mathbb{R}^d\rightarrow\mathbb{R}$,
\[
\nu^{(n)}_t(\omega,f)-\nu^{(n)}_0(\omega,f)
\stackrel{n\rightarrow+\infty}{\longrightarrow}0\qquad\mbox{a.s.};
\]

\textup{(ii)} if, moreover, a.s., every weak limit
$\nu^{(\infty)}(\omega,\mathrm{d}\alpha)$ of
$(\nu^{(n)}(\omega,\mathrm{d}\alpha))_{n\ge1}$
is the distribution of a Markov process with semigroup
$(Q_t^\omega)_{t\ge0}$, then, a.s., $\nu^{(\infty)}(\omega,\mathrm{d}\alpha)$
is the distribution of a stationary process.
\end{lemme}
\begin{pf}
(i) Let $f\dvtx\mathbb{R}^d\rightarrow\mathbb{R}$ be a bounded continuous
function. Since $\bar{X}_t^{(k)}=\bar{X}_{\Gamma_{N(k,t)}}$, we
have
\[
{\nu^{(n)}_t}(\omega,f)-{\nu^{(n)}_0}(\omega,f)=\frac{1}{H_n}
\sum_{k=1}^n \eta_k
\bigl(f\bigl(\bar{X}_{\Gamma_{N(k-1,t)}}\bigr)-f(\bar{X}_{\Gamma_{k-1}}) \bigr).
\]
From the very definition of $N(n,T)$ and $\tau(n,T)$, one checks
that $N(k-1,T)\le n-1$ if and only if $\tau(n,T)\ge k$. Then,
\begin{eqnarray*}
\frac{1}{H_n}\sum_{k=1}^n \eta_k f\bigl(\bar{X}_{\Gamma_{k-1}}\bigr)
&=& \frac{1}{H_n}\sum_{k=1}^{\tau(n,t)}\eta_{N(k-1,t)+1}
f\bigl(\bar{X}_{\Gamma_{N(k-1,t)}}\bigr)
\\
&&{} +\frac{1}{H_n}\sum_{k=1}^n\eta_k f(\bar{X}_{\Gamma_{k-1}})
1_{\{k-1\notin{N}(\{0,\ldots,n\},t)\}}.
\end{eqnarray*}
It follows that
\begin{eqnarray*}
{\nu^{(n)}_t}(\omega,f)-{\nu^{(n)}_0}(\omega,f)
&=& \frac{1}{H_n}\sum_{k=1}^{\tau(n,t)}
\bigl(\eta_k-\eta_{{N}(k-1,t)+1}\bigr)
f\bigl(\bar{X}_{\Gamma_{N(k-1,t)}}\bigr)
\\
&&{} +\frac{1}{H_n}\sum_{\tau(n,t)+1}^n\eta_k
f\bigl(\bar{X}_{\Gamma_{N(k-1,t)}}\bigr)
\\
&&{} - \frac{1}{H_n}\sum_{k=1}^n\eta_k
f(\bar{X}_{\Gamma_{k-1}})1_{\{k-1\notin{N}(\{0,\ldots,n\},t)\}}.
\end{eqnarray*}
Then, since $f$ is bounded and since
\begin{eqnarray*}
\sum_{k=1}^n\eta_k
1_{\{k-1\notin{N}(\{0,\ldots,n\},t)\}}
&=& \sum_{k=1}^n\eta_k-\sum_{k=1}^{\tau(n,t)}\eta_{{N}(k-1,t)+1}
\\
&\le & \sum_{k=1}^{\tau(n,t)}
\bigl|\eta_k-\eta_{{N}(k-1,t)+1}\bigr|
+ \sum_{k=\tau(n,t)+1}^n\eta_k,
\end{eqnarray*}
we deduce that
\[
\bigl|{\nu^{(n)}_t}(\omega,f)-{\nu^{(n)}_0}(\omega,f) \bigr|
\le 2\|f\|_\infty\Biggl(\frac{1}{H_n}\sum_{k=1}^{\tau(n,t)}
\bigl|\eta_k-\eta_{{N}(k-1,t)+1}\bigr|
+ \frac{1}{H_n}\sum_{k=\tau(n,t)+1}^n\eta_k \Biggr).
\]
Hence, we have to show that the sequences of the right-hand side
of the preceding inequality tend to 0. On the one hand, we observe
that
\[
\bigl|\eta_k-\eta_{{N}(k-1,t)+1}\bigr|
\le \sum_{\ell=k+1}^{N(k-1,T)+1}|\eta_\ell
-\eta_{\ell-1}|\le\max_{\ell\ge k+1}
\frac{|\Delta{\eta_\ell}|}{\gamma_\ell}
\sum_{\ell =k}^{N(k-1,T)+1}\gamma_\ell.
\]
Using the fact that $\sum_{\ell=k}^{N(k-1,T)+1}\gamma_\ell\le T+\gamma_1$
and condition (\ref{invarianceautomatique}) yields
\[
\frac{1}{H_n}\sum_{k=1}^{\tau(n,t)}
\bigl|\eta_k-\eta_{{N}(k-1,t)+1}\bigr|
\stackrel{n\rightarrow+\infty}{\longrightarrow}0.
\]
On the other hand, by $\mathbf{(C_{3,\varepsilon})}$, we have
\[
\frac{1}{H_n}\sum_{k=\tau(n,T)+1}^n\eta_k
\le
\frac{C}{H_n^{1-\varepsilon}}\sum_{k=\tau(n,T)+1}^{n}\gamma_k\le
\frac{CT}{H_n^{1-\varepsilon}}
\stackrel{n\rightarrow+\infty}{\longrightarrow}0
\qquad\mbox{a.s.},
\]
which completes the proof of (i).

(ii) Let $\mathbb{Q}_+$ denote the set of non-negative
rational numbers. Let $(f_\ell)_{\ell\ge1}$ be an everywhere dense
sequence in $\mathcal{C}_K(\mathbb{R}^d)$ endowed with the topology of
uniform convergence on compact sets. Since $\mathbb{Q}_+$ and
$(f_\ell)_{\ell\ge1}$ are countable, we derive from (i) that
there exists $\tilde{\Omega}\subset\Omega$ such that
\mbox{$\mathbb{P}(\tilde{\Omega})=1$} and such that for every
$\omega\in\tilde{\Omega}$, every $t\in\mathbb{Q}_+$ and
every $\ell\ge1$,
\[
\nu^{(n)}_t(\omega,f_\ell)-\nu^{(n)}_0(\omega,f_\ell)
\stackrel{n\rightarrow+\infty}{\longrightarrow}0.
\]
Let $\omega\in\tilde{\Omega}$ and let
$\nu^{(\infty)}(\omega,\mathrm{d}\alpha)$ denote a weak limit of
$(\nu^{(n)}(\omega,\mathrm{d}\alpha))_{n\ge1}$. We have
\[
\nu^{(\infty)}_t(\omega,f_\ell)=\nu^{(\infty)}_0
(\omega,f_\ell) \qquad\forall t\in\mathbb{Q}_+ \ \forall\ell\ge1
\]
and we easily deduce
that
\[
\nu^{(\infty)}_t(\omega,f)=\nu^{(\infty)}_0(\omega,f)\qquad
\forall t\in\mathbb{R}_+\ \forall f\in\mathcal{C}_K(\mathbb{R}^d).
\]
Hence, if
$\nu^{(\infty)}(\omega,\mathrm{d}\alpha)$ is the distribution of a Markov
process $(Y_t)$ with semigroup $(Q_t^\omega)_{t\ge0}$, we have,
for all $f\in\mathcal{C}_K(\mathbb{R}^d)$,
\[
\int Q_t^\omega f(x)\nu_0^{(\infty)}(\omega,\mathrm{d}x)
= \int f(x)\nu_0^{(\infty)}(\omega,\mathrm{d}x) \qquad\forall t\ge0.
\]
$\nu_0^{(\infty)}(\omega,\mathrm{d}x)$ is then an
invariant distribution for $(Y_t)$. This completes the proof.
\end{pf}

\subsection[Proof of Theorem 1]{Proof of Theorem \protect\ref{prinfonc1}}

Thanks to Lemma \ref{caracterisation}(a) applied with
$E=\mathbb{D}(\mathbb{R}_+,\mathbb{R}^d)$ and $d$ defined by
(\ref{distancecomplete}),
%
\begin{equation}\label{suffis}
\nu^{(n)}(\omega,\mathrm{d}\alpha)\stackrel{(Sk)}
{\Longrightarrow}\mathbb{P}_{\nu_0}(\mathrm{d}\alpha)
\qquad \mbox{a.s.}\Longleftrightarrow{{\nu}^{(n)}}(\omega,F)
\stackrel{n\rightarrow+\infty}{\longrightarrow}
\int F(x)\mathbb{P}_{\nu_0}(\mathrm{d}x)\qquad \mbox{a.s.}
\end{equation}
\noindent for every bounded Lipschitz functional
$F\dvtx\mathbb{D}(\mathbb{R}_+,\mathbb{R}^d)\rightarrow\mathbb{R}$.
Now, consider such a
functional. By the assumptions of Theorem \ref{prinfonc1}, we know
that a.s., $(\nu^{(n)}_0(\omega,\mathrm{d}x))_{n\ge1}$ converges weakly
to ${\nu_0}$. Set $\Phi^F(x):=\mathbb{E}[F(X^x)]$, $x\in\mathbb{R}^d$. By
Lemma \ref{lemmeflow}(a), $\Phi^F$ is a bounded continuous
function on $\mathbb{R}^d$. It then
follows from $\mathbf{(C_{0,1})}$ that
\[
\frac{1}{H_n}\sum_{k=1}^n\eta_k\Phi^F
\bigl(\bar{X}_0^{(k-1)}\bigr)
\stackrel{n\rightarrow+\infty}{\longrightarrow}
\int\Phi^F(x){\nu_0}(\mathrm{d}x)
= \int F(x)\mathbb{P}_{\nu_0}(\mathrm{d}x)\qquad \mbox{a.s.}
\]
Hence, the right-hand side of (\ref{suffis}) holds for $F$ as soon
as
%
%
\begin{equation}\label{ucops3}
\frac{1}{H_n}\sum_{k=1}^n\eta_k\bigl (F\bigl(\bar{X}^{(k-1)}\bigr)
- \Phi^F\bigl(\bar{X}_0^{(k-1)}\bigr) \bigr)
\stackrel{n\rightarrow+\infty}{\longrightarrow}0\qquad \mbox{a.s.}
\end{equation}
Let
us prove (\ref{ucops3}). First, let $T>0$ and let
$F_T$ be defined by (\ref{functionalt}). By Lemma \ref{fundamental},
%
%
\begin{equation}\label{ucops2}
\frac{1}{H_n}\sum_{k=1}^n\eta_k
F_T\bigl(\bar{X}^{(k-1)}\bigr)
- \frac{1}{H_n}\sum_{k=1}^n\eta_k\mathbb{E}
\bigl[F_T\bigl(\bar{X}^{(k-1)}\bigr)/\bar{\mathcal{F}}_{\Gamma_{k-1}}\bigr]
\stackrel{n\rightarrow+\infty}{\longrightarrow}0\qquad \mbox{a.s.}
\end{equation}
With the notation of Lemma \ref{lemmeflow}(b), we derive from
assumption $\mathbf{(C_2)}$(i) that
\[
\mathbb{E}\bigl[F_T\bigl(\bar{X}^{(k-1)}\bigr)/
\bar{\mathcal{F}}_{\Gamma_{k-1}}\bigr]
=\Phi^{F_T}_k\bigl(\bar{X}_0^{(k-1)}\bigr).
\]
Let $N\in\mathbb{N}$. On one hand, by Lemma \ref{lemmeflow}(b),
%
%
\begin{equation}\label{ucops}
\frac{1}{H_n}\sum_{k=1}^n\eta_k
\bigl(\Phi^{F_T}_k\bigl(\bar{X}_0^{(k-1)}\bigr)
- \Phi^{F_T}\bigl(\bar{X}_0^{(k-1)}\bigr) \bigr)
1_{\{|\bar{X}_0^{(k-1)}|\le N\}}
\stackrel{n\rightarrow+\infty}{\longrightarrow}0\qquad \mbox{a.s.}
\end{equation}
On
the other hand, the tightness of
$(\nu^{(n)}_0(\omega,\mathrm{d}x))_{n\ge1}$ on $\mathbb{R}^d$ yields
\[
\psi(\omega,N):=\sup_{n\ge1}\bigl(\nu^{(n)}_0(\omega
,(B(0,N)^c))\bigr)\stackrel{N\rightarrow+\infty}{\longrightarrow}0\qquad \mbox{a.s.}
\]
It follows that, a.s.,
%
%
\begin{eqnarray}\label{fonc622}
&& \sup_{n\ge 1} \Biggl(\frac{1}{H_n}\sum_{k=1}^n\eta_k
\bigl|\Phi^{F_T}_k\bigl(\bar{X}_0^{(k-1)}\bigr)
- \Phi^{F_T}\bigl(\bar{X}_0^{(k-1)}\bigr) \bigr|
1_{\{|\bar{X}_0^{(k-1)}|> N\}} \Biggr)
\nonumber\\[-8pt]
\\[-8pt]
&&\quad \le 2\|F\|_\infty\psi(\omega,N)
\stackrel{N\rightarrow+\infty}{\longrightarrow}0.
\nonumber
\end{eqnarray}
Hence, a
combination of (\ref{ucops}) and (\ref{fonc622}) yields
%
%
\begin{equation}\label{numero10}
\forall T>0\qquad
\frac{1}{H_n}\sum_{k=1}^n\eta_k
\bigl(\Phi^{F_T}_k\bigl(\bar{X}_0^{(k-1)}\bigr)
- \Phi^{F_T}\bigl(\bar{X}_0^{(k-1)}\bigr) \bigr)
\stackrel{n\rightarrow+\infty}{\longrightarrow}0\qquad \mbox{a.s.}
\end{equation}
Finally, let $(T_\ell)_{\ell\ge1}$ be a sequence of positive
numbers such that, $T_\ell\rightarrow+\infty$ when
$\ell\rightarrow+\infty$. Combining (\ref{numero10}) and
(\ref{ucops2}), we obtain that, a.s., for every $\ell\ge1$,
\begin{eqnarray*}
&& \limsup_{n\rightarrow+\infty}
\Biggl|\frac{1}{H_n}\sum_{k=1}^n\eta_k
\bigl(F\bigl(\bar{X}^{(k-1)}\bigr)
- \Phi^F\bigl(\bar{X}^{(k-1)}\bigr) \bigr) \Biggr|
\\
&&\quad \le \limsup_{n\rightarrow+\infty}
\Biggl|\frac{1}{H_n}\sum_{k=1}^n\eta_k
\bigl(F\bigl(\bar{X}^{(k-1)}\bigr)
- F_{T_\ell}\bigl(\bar{X}^{(k-1)}\bigr) \bigr) \Biggr|
\\
&&\qquad {} +\limsup_{n\rightarrow+\infty}
\Biggl|\frac{1}{H_n}\sum_{k=1}^n\eta_k
\bigl(\Phi^{F_{T_\ell}}\bigl(\bar{X}_0^{(k-1)}\bigr)
- \Phi^{F}\bigl(\bar{X}_0^{(k-1)}\bigr) \bigr) \Biggr|.
\end{eqnarray*}
By the definition of ${d}$,
$|F-F_{T_\ell}|\le\mathrm{e}^{-T_\ell}$. Then, a.s.,
\[
\limsup_{n\rightarrow+\infty}
\Biggl|\frac{1}{H_n}\sum_{k=1}^n\eta_k
\bigl(F\bigl(\bar{X}^{(k-1)}\bigr)
- \Phi^F\bigl(\bar{X}_0^{(k-1)}\bigr) \bigr) \Biggr|
\le 2\mathrm{e}^{-T_\ell}\qquad\forall\ell\ge1.
\]
Letting $\ell\rightarrow+\infty$ implies (\ref{ucops3}).

The generalization to non-bounded functionals in Theorem
\ref{prinfonc1} is then derived from (\ref{equiintfuncassump}) and
from a uniform integrability argument.

\subsection[Proof of Theorem 2]{Proof of Theorem \protect\ref{prinfonc2}}

(i) We want to prove that the conditions of Lemma
\ref{caracterisation}(b) are fulfilled. Since
$(\nu^{(n)}_0(\omega,\mathrm{d}x))_{n\ge1}$ is supposed to be a.s. tight,
one can check that for every bounded Lipschitz functional
$F\dvtx\mathbb{D}(\mathbb{R}_+,\mathbb{R}^d)\rightarrow\mathbb{R}$,
(\ref{ucops3}) is still
valid. Then, let $(F_{\ell})_{\ell\ge1}$ be a sequence of bounded
Lipschitz functionals. There exists $\tilde{\Omega}\subset\Omega$
with $\mathbb{P}(\tilde{\Omega})=1$ such that for every
$\omega\in\tilde{\Omega}$, $(\nu^{(n)}_0(\omega,\mathrm{d}x))_{n\ge1}$ is
tight and
%
\begin{equation}\label{ucopsfunc4}
\frac{1}{H_n}\sum_{k=1}^n\eta_k
\bigl(F_{\ell}\bigl(\bar{X}^{(k-1)}(\omega)\bigr)
- \Phi^{F_{\ell}}\bigl(\bar{X}_0^{(k-1)}(\omega)\bigr) \bigr)
\stackrel{n\rightarrow+\infty}{\longrightarrow}0\qquad\forall\ell\ge1.
\end{equation}
Let $\omega\in\tilde{\Omega}$ and let
$\phi_\omega\dvtx\mathbb{N}\mapsto\mathbb{N}$ be an increasing function. As
$(\nu^{(\phi_\omega(n))}_0(\omega,\mathrm{d}x))_{n\ge1}$ is
tight, there exists a convergent subsequence
$(\nu^{(\phi_\omega\circ\psi_\omega(n))}_0(\omega,\mathrm{d}x))_{n\ge1}$.
We denote its weak limit by $\nu_\infty$. Since $\Phi^{F_{\ell}}$
is continuous for every $\ell\ge1$ (see Lemma \ref{lemmeflow}(a)),
\[
\nu^{(\phi_\omega\circ\psi_\omega(n))}_0(\omega,\Phi^{F_{\ell}})
\stackrel{n\rightarrow+\infty}{\longrightarrow}
\nu_\infty(\Phi^{F_{\ell}}) =\int F_{\ell}(\alpha)
\mathbb{P}_{\nu_\infty}(\mathrm{d}\alpha)\qquad \forall\ell\ge1.
\]
We then derive from (\ref{ucopsfunc4}) that
for every $\ell\ge1$
\[
\nu^{(\phi_\omega\circ\psi_\omega(n))}
(\omega,F_{\ell})\stackrel{n\rightarrow+\infty}{\longrightarrow}
\int F_{\ell}(\alpha)\mathbb{P}_{\nu_\infty}(\mathrm{d}\alpha).
\]
It follows that the conditions of Lemma
\ref{caracterisation}(b) are fulfilled with
$\mathcal{U}=\{\mathbb{P}_\mu,\mu\in\mathcal{I}\}$, where
\[
\mathcal{I}= \biggl\{\mu\in\mathcal{P}(\mathbb{R}^d),
\exists\omega\in\tilde{\Omega}
\mbox{ and an increasing function }
\phi\dvtx \mathbb{N}\mapsto\mathbb{N},
\mu=\lim_{n\rightarrow+\infty}\nu^{(\phi(n))}
(\omega,\mathrm{d}\alpha)\biggr\}.
\]
Hence, by Lemma \ref{caracterisation}(b), we deduce that
$(\nu^{(n)}(\omega,\mathrm{d}\alpha))_{n\ge1}$ is a.s. tight with
$\mathcal{U}$-valued limits.

Finally, Theorem \ref{prinfonc2}(ii) is a consequence of condition
(\ref{invarianceautomatique}) and Lemma \ref{funclemme3}(ii).

\section{Path-dependent option pricing in stationary stochastic
volatility models}\label{foncfinance}

In this section, we propose a simple and efficient method to price
options in stationary stochastic volatility (SSV) models. In most
stochastic volatility (SV) models, the volatility is a
mean reverting process. These processes are generally ergodic
with a unique invariant distribution
(the Heston model or the BNS model for instance
(see below) but also the SABR model (see Hagan \textit{et al.}
\cite{hagan})$,\ldots)$.
However, they are usually considered
in SV models under a non-stationary regime, starting from a
deterministic value (which usually turns out to be the mean of
their invariant distribution). However, the instantaneous
volatility is not easy to observe on the market since it is not a
traded asset. Hence,
it seems to be more natural to assume that it evolves under its
stationary regime than to give it a deterministic value at time~0.
\footnote{When one has sufficiently close observations of the stock
price, it is in fact possible to derive a rough idea of the size of the
volatility from the variations of the stock price (see, e.g., Jacod
\cite{jacodpower}). Then, using this information, a good compromise
between a deterministic initial value and the stationary case may be to assume
that the distribution $\mu_0$ of the volatility at time 0 is
concentrated around the estimated value (see Section \ref{remarqueabs}
for application of our algorithm in this case).}

From a purely calibration viewpoint, considering an SV model in
its SSV regime will not modify the set of parameters used to
generate the implied volatility surface, although it will modify
its shape, mainly for short maturities. This effect can in fact be
an asset of the SSV approach since it may correct some observed
drawbacks of some models (see, e.g., the Heston model below).

From a numerical point of view, considering SSV models is no
longer an obstacle, especially when considering multi-asset
models (in the unidimensional case, the stationary distribution
can be made more or less explicit like in the Heston model; see
below) since our algorithm is precisely devised to compute by
simulation some expectations of functionals of processes under
their stationary regime, even if this stationary regime cannot be
directly simulated.

As a first illustration (and a benchmark) of the method, we will
describe in detail the algorithm for the pricing of Asian options
in a Heston model. We will then show in our numerical results to
what extent it differs, in terms of smile and skew,
from the usual SV Heston model for short maturities. Finally, we will
complete this section with a numerical test on Asian options in the BNS
model where the
volatility is driven by a tempered stable subordinator. Let us
also mention that this method can be applied to other fields of
finance like interest rates, and commodities and energy
derivatives where mean-reverting processes play an important role.

\subsection{Option pricing in the Heston SSV model}\label{methodpricing}

We consider a Heston stochastic
volatility model. The dynamic of the asset price process
$(S_t)_{t\ge0}$ is given by $S_0=s_0$ and
\begin{eqnarray*}
\mathrm{d}S_t &=& S_t\bigl(r\,\mathrm{d}t
+ \sqrt{(1-\rho^2)v_t}\,\mathrm{d}W^1_t
+ \rho\sqrt{v_t} \,\mathrm{d}W_t^2\bigr),
\\
\mathrm{d}v_t &=& k({\theta}-v_t)\,\mathrm{d}t
+ \varsigma\sqrt{v_t}\,\mathrm{d}W_t^2,
\end{eqnarray*}
where $r$ denotes the interest rate, $(W^1,W^2)$ is a standard
two-dimensional Brownian motion, $\rho\in[-1,1]$ and $k$, $\theta$
and $\varsigma$ are some non-negative numbers. This model was
introduced by Heston in 1993 (see Heston \cite{heston}). The equation for
$(v_t)$ has a unique (strong) pathwise continuous solution living
in $\mathbb{R}_+$. If, moreover, $2k\theta>\varsigma^2$, then $(v_t)$
is a positive process (see Lamberton and Lapeyre \cite{lamblapeyre}). In
this case, $(v_t)$ has a unique invariant probability ${\nu_0}$. Moreover,
${\nu_0}=\gamma(a,b)$ with $a=(2k)/\varsigma^2$ and
$b=(2k\theta)/\varsigma^2$.
In the following, we will assume that
$(v_t)$ is in its stationary regime, that is, that
\[
\mathcal{L}(v_0)={\nu_0}.
\]

\subsubsection{Option price and stationary processes}

Using our procedure to price options in this model naturally needs
to express the option price as the expectation of a
functional of a stationary stochastic process.

\textbf{Na\"{\i}ve method.}
(may work) Since $(v_t)_{t\ge0}$ is stationary, the
first idea is to express the option price as
the expectation of a functional of
$(v_t)_{t\ge0}$: by It\^{o} calculus, we have
%
\begin{equation}\label{hestoneq}
S_t=s_0\exp\biggl(\biggl(rt-\tfrac{1}{2}\int_0^tv_s\,\mathrm{d}s\biggr)
+\rho\int_0^t\sqrt{v_s}\,\mathrm{d}W_s^2
+ \sqrt{1-\rho^2}\int_0^t\sqrt{v_s}\,\mathrm{d}W^1_s \biggr).
\end{equation}
Since
\[
\int_0^t\sqrt{v_s}\,\mathrm{d}W_s^2
= \Lambda(t,(v_t)):=\frac{v_t-v_0-k\theta t+k\int_0^t v_s\,\mathrm{d}s}{\varsigma},
\]
it follows by setting $M_t=\int_0^t\sqrt{v_s}\,\mathrm{d}W^1_s$ that
%
%
\begin{equation}\label{Stexpress}
S_t=\Psi(t,(v_s),(M_s)),
\end{equation}
where $\Psi$ is given for
every $t\ge0$, $u$ and $w\in\mathcal{C}(\mathbb{R}_+,\mathbb{R})$ by
\[
\Psi(t,u,w)=s_0\exp\biggl(\biggl(rt-\tfrac{1}{2}\int_0^t
u(s)\,\mathrm{d}s\biggr)+\rho\Lambda(t,u)+\sqrt{1-\rho^2}w(t) \biggr).
\]
Then, let
$F\dvtx\mathcal{C}(\mathbb{R}_+,\mathbb{R})\rightarrow\mathbb{R}$ be a
non-negative measurable
functional. Conditioning by $\mathcal{F}^{W^2}_T$ yields
\[
\mathbb{E}[F_T((S_t)_{t\ge0})]=\mathbb{E}[\tilde{F}_T((v_t)_{t\ge0})],
\]
where, for every $u\in\mathcal{C}(\mathbb{R}_+,\mathbb{R})$,
\[
\tilde{F}_T(u)=\mathbb{E}\biggl[F_T\biggl(\biggl(\Psi\biggl(t,u,
\int_0^t\sqrt{u(s)}\,\mathrm{d}W_s^1\biggr)\biggr)_{t\ge0}\biggr)\biggr].
\]
For some particular options such as the European call or put (thanks
to the Black--Scholes formula), the functional $\tilde{F}$ is
explicit. In those cases, this method seems to be very efficient
(see Panloup \cite{panloup2} for numerical results). However, in the
general case, the computation of $\tilde{F}$ will need some
Monte Carlo methods at each step. This approach is then very
time-consuming in general -- that is why we are going
to introduce another representation of the option as a functional of a
stationary process.

\textbf{General method.}
(always works) We express the option premium as the
expectation of a functional of a two-dimensional stationary
stochastic process. This method is based on the following idea.
Even though $(v_t,M_t)$ is not stationary, $(S_t)$ can be expressed as a
functional of a stationary process $(v_t,y_t)$. Indeed, consider
the following SDE given by
%
%
\begin{equation}
\cases{%
\mathrm{d}y_t=-y_t\,\mathrm{d}t+\sqrt{v_t}\,\mathrm{d}W_t^1,\cr
\mathrm{d}v_t=k({\theta}-v_t)\,\mathrm{d}t+\varsigma\sqrt{v_t}\,\mathrm{d}W_t^2.}
\end{equation}
First, one checks that the SDE has a unique strong solution and that
assumption $\mathbf{(S_1)}$ is fulfilled with
$V(x_1,x_2)=1+x_1^2+x_2^2$. This
ensures the existence of an invariant distribution
$\tilde{{\nu}}_0$ for the SDE (see, e.g., Pag\`es \cite{bib14}). Then,
since $(v_t)$ is positive and has a unique
invariant distribution, the uniqueness of the invariant
distribution follows. Then, assume that
$\mathcal{L}(y_0,v_0)=\tilde{{\nu}}_0$. Since
$(v_t,M_t)=(v_t,y_t-y_0+\int_0^ty_s\,\mathrm{d}s)$, we have, for every
positive measurable functional
$F\dvtx\mathcal{C}(\mathbb{R}_+,\mathbb{R})\rightarrow\mathbb{R}$,
%
%
\begin{eqnarray}\label{equaltyfunc}
\mathbb{E}[F_T((S_t)_{t\ge0})]&=& \mathbb{E}[F_T((\psi(t,v_t,M_t))_{t\ge0})]
\nonumber\\[-8pt]
\\[-8pt]
&=& \mathbb{E}_{\tilde{{\nu_0}}}
\biggl[F_T\biggl(\biggl(\psi\biggl(t,v_t,y_t-y_0
+ \int_0^ty_s\,\mathrm{d}s\biggr)\biggr)_{t\ge0}\biggr)\biggr],
\nonumber
\end{eqnarray}
where $\mathbb{P}_{\tilde{{\nu_0}}}$ is the stationary distribution of
the process $(v_t,y_t)$. Every option
price can then be expressed as the expectation of an explicit
functional of a stationary process. We will develop this second general
approach in the numerical tests below.
\begin{Remarque}
The idea of the second method holds for every stochastic
volatility model for which~$(S_t)$ can be written as follows:
%
%
\begin{equation}\label{repoptionprice}
S_t=\Phi\Biggl(t,v_t,\sum_{i=1}^p\int_0^t
h_i(|v_s|)\,\mathrm{d}Y^i_s\Biggr),
\end{equation}
where, for every $i\in\{1,\ldots,p\}$,
$h_i\dvtx\mathbb{R}_+\rightarrow\mathbb{R}$ is a positive function such that
$h_i(x)=o(|x|)$ as $|x|\rightarrow+\infty$, $(Y_t^i)$ is a
square-integrable centered L\'{e}vy process and $(v_t)$ is a
mean reverting stochastic process solution to a L\'{e}vy driven SDE.

In some complex models, showing the uniqueness of the invariant
distribution may be difficult. In fact, it is important to note
at this stage that the uniqueness of the invariant distribution
for the couple $(v_t,y_t)$ is not required. Indeed, by
construction, the local martingale $(M_t)$ does not depend on the
choice of $y_0$. It follows that if
$\mathcal{L}(y_0,v_0)=\tilde{\mu}$,
with $\tilde{\mu}$ constructed such that
$\mathcal{L}(v_0)={\nu_0}$, (\ref{equaltyfunc}) still holds.
This implies that it is only necessary that uniqueness holds
for the invariant distribution of the stochastic volatility
process.
\end{Remarque}

\subsubsection{Numerical tests on Asian options}

We recall that $(v_t)$ is a Cox--Ingersoll--Ross process. For this
type of processes, it is well known that the genuine Euler scheme
cannot be implemented since it does not preserve the non-negativity of
the $(v_t)$.
That is why some specific discretization schemes
have been studied by several authors (Alfonsi \cite{alfonsi},
Deelstra and Delbaen \cite{delstra} and Berkaoui \textit{et al.}
\cite{berkaoui,diop}). In this paper, we consider the scheme
studied by the last authors in a decreasing step framework. We denote
it by
$(\bar{v}_t)$. We set $\bar{v}_0=x>0$ and
\[
\bar{v}_{\Gamma_{n+1}}= \bigl|\bar{v}_{\Gamma_n}+k\gamma_{n+1}
(\theta-\bar{v}_{\Gamma_n})+ \varsigma\sqrt{\bar{v}_{\Gamma_n}}
(W^2_{\Gamma_{n+1}}-W^2_{\Gamma_n}) \bigr|.
\]
We also introduce the stepwise constant Euler scheme $(\bar{y}_t)$ of
$(y_t)_{t\ge0}$ defined by
\[
\bar{y}_{\Gamma_{n+1}}=\bar{y}_{{\Gamma_n}}-\gamma_{n+1}
\bar{y}_{\Gamma_n}+\sqrt{\bar{v}_{\Gamma_n}}
(\tilde{W}^1_{\Gamma_{n+1}}-\tilde{W}^1_{\Gamma_{n}}),\qquad
\bar{y}_0=y\in\mathbb{R}^d.
\]
Denote by $(\bar{v}_t^{(k)})$ and $(\bar{y}_t^{(k)})$ the shifted
processes defined by $\bar{v}_t^{(k)}:=\bar{v}_{\Gamma_k+t}$ and
$\bar{y}_t^{(k)}=\bar{y}_{\Gamma_k+t}$, and let
$(\nu^{(n)}(\omega,\mathrm{d}\alpha))_{n\ge1}$ be the sequence of empirical
measures defined by
\[
\nu^{(n)}(\omega,\mathrm{d}\alpha)
= \frac{1}{H_n}\sum_{k=1}^n\eta_k
1_{\{(\bar{v}^{(k-1)},\bar{y}^{(k-1)})\in d\alpha\}}.
\]
The specificity of both the model and the Euler scheme implies that
Theorems \ref{prinfonc1} and \ref{prinfonc2} cannot be directly
applied here. However, a specific study using the fact that (\ref{djioe3})
holds for every compact set of $\mathbb{R}_+^*\times\mathbb{R}$ when
$2k\theta/\varsigma^2>1+2\sqrt{6}/\varsigma$ (see Theorem 2.2 of
Berkaoui \textit{et al.}
\cite{berkaoui} and Remark~\ref{djioe3}) shows that
\[
\nu^{(n)}(\omega,\mathrm{d}\alpha)\stackrel{n\rightarrow+\infty}
{\Longrightarrow}\mathbb{P}_{\tilde{\nu}_0}
(\mathrm{d}\alpha)\qquad \mbox{a.s.}
\]
when $2k\theta/\varsigma^2>1+2\sqrt{6}/\varsigma$. Details are left
to the reader.

Let us now state our numerical results obtained for the
pricing of
Asian options with this discretization. We denote by
$C_{as}(\nu_0,K,T)$ and $P_{as}(\nu_0,K,T)$ the Asian call and put
prices in the SSV Heston model. We have
\begin{eqnarray*}
C_{as}(\nu_0,K,T)=\mathrm{e}^{-rT}\mathbb{E}_{\nu_0}
\biggl[\biggl(\frac{1}{T}\int_0^TS_s\,\mathrm{d}s-K\biggr)_+\biggr]
\end{eqnarray*}
and
\[
P_{as}(\nu_0,K,T)=\mathrm{e}^{-rT}\mathbb{E}_{\nu_0}
\biggl[\biggl(K-\frac{1}{T}\int_0^TS_s\,\mathrm{d}s\biggr)_+\biggr].
\]
With the notation of (\ref{equaltyfunc}), approximating
$C_{as}(\nu_0,K,T)$ and $P_{as}(\nu_0,K,T)$ by our procedure needs
to simulate the sequences $(C_{as}^n)_{n\ge1}$ and
$(P_{as}^n)_{n\ge1}$ defined by
\begin{eqnarray*}
C_{as}^n &=& \frac{1}{H_n}\sum_{k=1}^n\eta_k\mathrm{e}^{-rT}
\biggl(\frac{1}{T}\int_0^T\Psi
\bigl(s,\bar{v}^{(k-1)},\bar{M}^{(k-1)}\bigr)\,\mathrm{d}s
-K \biggr)_+,
\\
P_{as}^n &=& \frac{1}{H_n}\sum_{k=1}^n\eta_k\mathrm{e}^{-rT}
\biggl(K-\frac{1}{T}\int_0^T\Psi
\bigl(s,\bar{v}^{(k-1)},\bar{M}^{(k-1)}\bigr)\,\mathrm{d}s \biggr)_+.
\end{eqnarray*}
These sequences can be computed by the method developed in
Section \ref{simfuncpropos}. Note that the specific properties
of the exponential function and the linearity of the integral
imply that $(\int_0^T\Psi(t,\bar{v}^{(n-1)},\bar{M}^{(n-1)})\,\mathrm{d}s)$
can be computed quasi-recursively.

\begin{table}[b]
\caption{Approximation of the Asian call price \label{callasian}}
\begin{tabular*}{\textwidth}{@{\extracolsep{\fill}} lccccccc@{}}
\hline
$K$ & 44 &45& 46&47 & 48&49&50 \\
\hline
Asian call (ref.) & \textbf{6.92} &\textbf{5.97}& \textbf{5.04}& \textbf{4.12} &\textbf{3.25} & \textbf{2.46} & \textbf{1.78}\\
$N=5\cdot10^4$& 6.89& 6.07& 5.07& 4.13& 3.18& 2.49& 1.77\\
$N=5\cdot10^5$& 6.90 & 6.02& 5.00& 4.11& 3.24& 2.46& 1.79\\
$N=5\cdot10^4$ (CP parity)& 6.92& 5.96& 5.04& 4.13& 3.26& 2.46& 1.78\\
$N=5\cdot10^5$ (CP parity)&6.92& 5.97& 5.04& 4.12& 3.25& 2.47& 1.78 
\end{tabular*}
\begin{tabular*}{\textwidth}{@{\extracolsep{\fill}} lcccccc@{}}
\hline
$K$ & 51 &52& 53&54 & 55&56 \\
\hline
Asian call (ref.)& \textbf{1.23} & \textbf{0.82} & \textbf{0.53}& \textbf{0.33} &\textbf{0.21}& \textbf{0.12}\\
$N=5\cdot10^4$&1.21& 0.81 & 0.51& 0.34& 0.22& 0.11\\
$N=5\cdot10^5$& 1.23& 0.82& 0.53& 0.33& 0.21& 0.13 \\
$N=5\cdot10^4$ (CP parity)& 1.23& 0.82& 0.53& 0.31& 0.21& 0.12\\
$N=5\cdot10^5$(CP parity)& 1.23& 0.82& 0.53& 0.33& 0.21 & 0.13\\
\hline
\end{tabular*}
\end{table}

Let us state our numerical results for the Asian call with
parameters
%
\begin{eqnarray}\label{paramfunc}
s_0 &=& 50,\qquad
r=0.05 ,\qquad
T=1 ,\qquad
\rho=0.5,\qquad
\nonumber\\[-8pt]
\\[-8pt]
\theta &=&0.01,\qquad
\varsigma=0.1,\qquad
k=2.
\nonumber
\end{eqnarray}
We
also assume that $K\in\{44,\ldots,56\}$ and choose the following
steps and weights: $\gamma_n=\eta_n=n^{-{1}/{3}}$. In Table
\ref{callasian}, we first state the reference value for the Asian
call price obtained for $N=10^8$ iterations. In the two following
lines, we state our results for $N=5.10^4$ and $N=5.10^5$
iterations. Then, in the last lines, we present the numerical
results obtained using the call-put parity
%
%
\begin{equation}
C_{as}(\nu_0,K,T)-P_{as}(\nu_0,S_0,K,T)=\frac{s_0}{rT}
(1- \mathrm{e}^{-rT}) -K\mathrm{e}^{-rT}
\end{equation}
as a means of variance reduction. The computation times for
$N=5.10^4$ and $N=5.10^5$ (using MATLAB with a Xeon 2.4
GHz processor) are about $5$ s and 51~s, respectively.
In particular, the complexity is
quasi-linear and the additional computations needed when we use
the call-put parity are negligible.

\subsection{Implied volatility surfaces of Heston SSV and SV
models}

Given a particular pricing model (with initial value $s_0$
and interest rate $r$) and its associated European call prices
denoted by $C_{\mathrm{eur}}(K,T)$, we recall that the implied
volatility surface is the graph of the function
$(K,T)\mapsto\sigma_{\mathrm{imp}}(K,T)$, where
$\sigma_{\mathrm{imp}}(K,T)$ is defined for every maturity $T>0$ and
strike $K$ as the unique solution of
\[
C_{BS}(s_0,K,T,r,\sigma_{\mathrm{imp}}(K,T))=C_{\mathrm{eur}}(K,T),
\]
where $C_{BS}(s_0,K,T,r,\sigma)$ is the price of the European call
in the Black--Scholes model with parameters $s_0$, $r$ and
$\sigma$. When $C_{\mathrm{eur}}(K,T)$ is known, the value of
$\sigma_{\mathrm{imp}}(K,T)$ can be numerically computed
using the Newton method or by dichotomy if the first method is not
convergent.

In this last part, we compare the implied volatility surfaces
induced by the SSV and SV Heston models where we suppose that the
initial value of $(v_t)$ in the SV Heston model is the mean of the
invariant distribution, that is, we suppose that
$v_0=\theta$.
\footnote{This choice is the most usual in practice.}
We also assume that the parameters are those of (\ref{paramfunc}),
except the correlation coefficient $\rho$.

\begin{figure}

\includegraphics{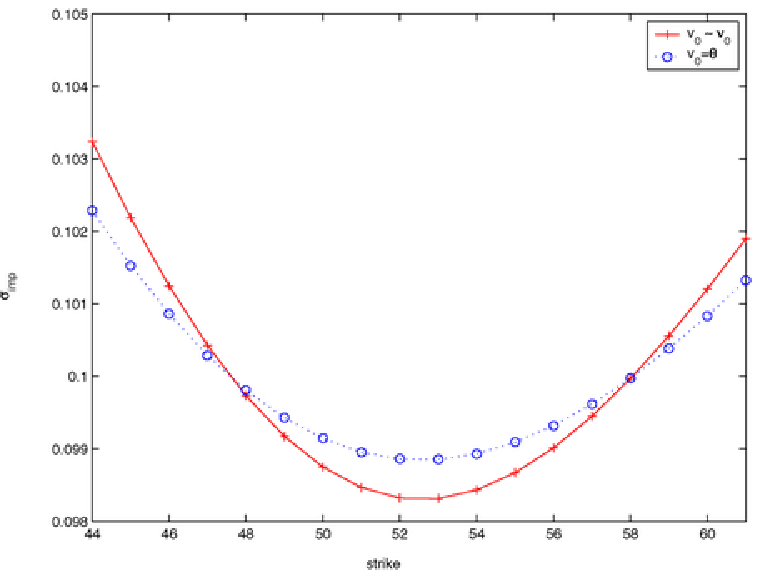}

\caption{$\rho=0$, $K\mapsto\sigma_{\mathrm{imp}}(K,1)$.}
\label{foncfigure1}
\end{figure}

\begin{figure}[b]

\includegraphics{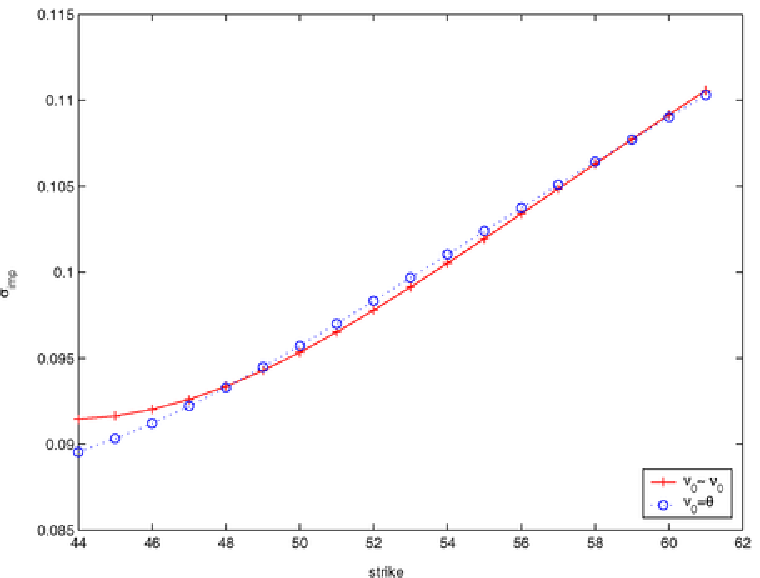}

\caption{$\rho=0.5$, $K\mapsto\sigma_{\mathrm{imp}}(K,1)$.}
\label{foncfigure2}
\end{figure}

\begin{figure}

\includegraphics{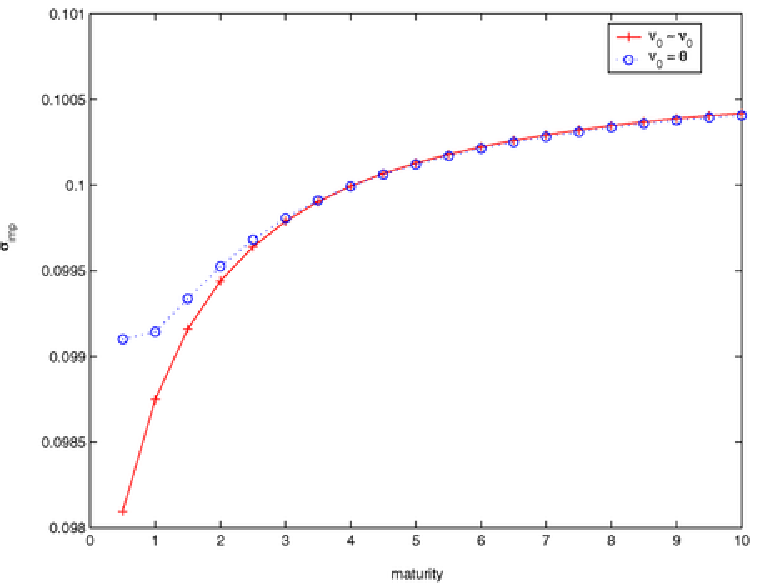}

\caption{$\rho=0$, $T\mapsto\sigma_{\mathrm{imp}}(50,T)$.}
\label{foncfigure3}
\end{figure}

\begin{figure}[b]

\includegraphics{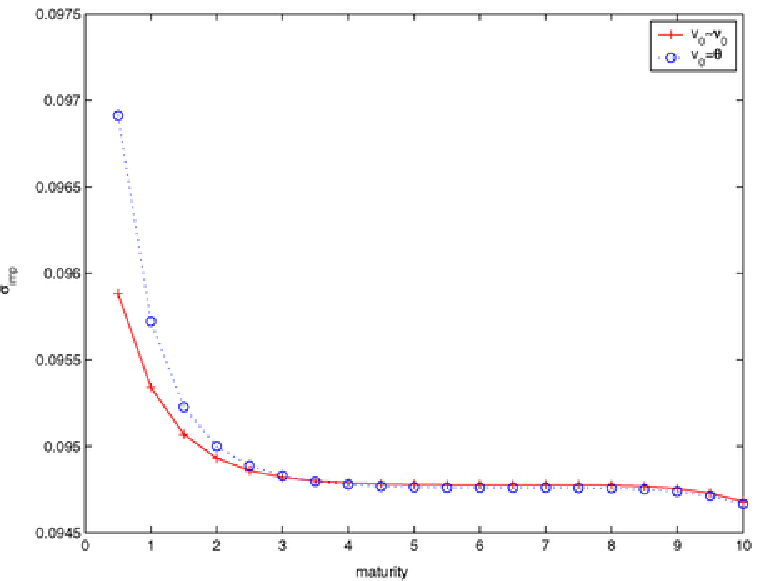}

\caption{$\rho=0.5$, $T\mapsto\sigma_{\mathrm{imp}}(50,T)$.}
\label{foncfigure4}
\end{figure}

In Figures \ref{foncfigure1} and \ref{foncfigure2},
the volatility curves obtained when $T=1$ are depicted, whereas in Figures
\ref{foncfigure3} and \ref{foncfigure4}, we set the strike $K$
at $K=50$ and let the time vary. These representations show that
when the maturity is long, the differences between the
SSV and SV Heston models vanish. This is a consequence of the
convergence of the stochastic volatility to its stationary regime
when $T\rightarrow+\infty$.

\begin{figure}

\includegraphics{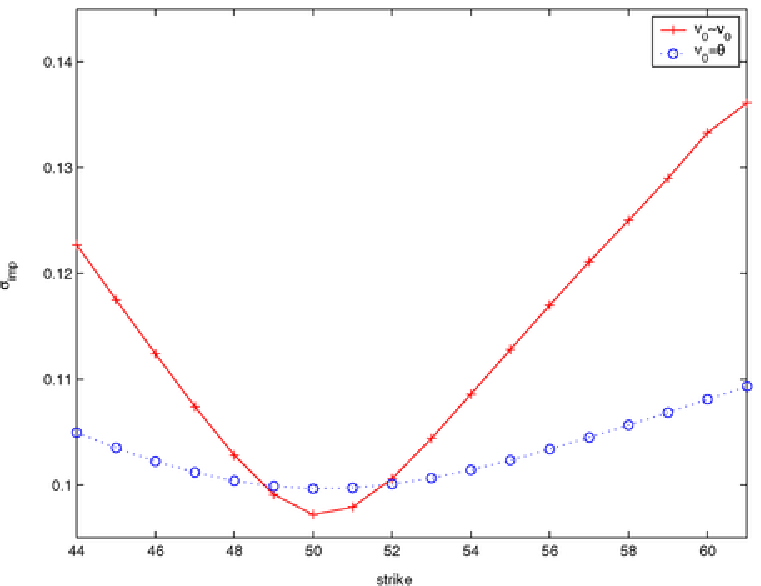}

\caption{$\rho=0$, $T\mapsto\sigma_{\mathrm{imp}}(50,T)$.}
\label{foncfigure5}
\end{figure}

\begin{figure}[b]

\includegraphics{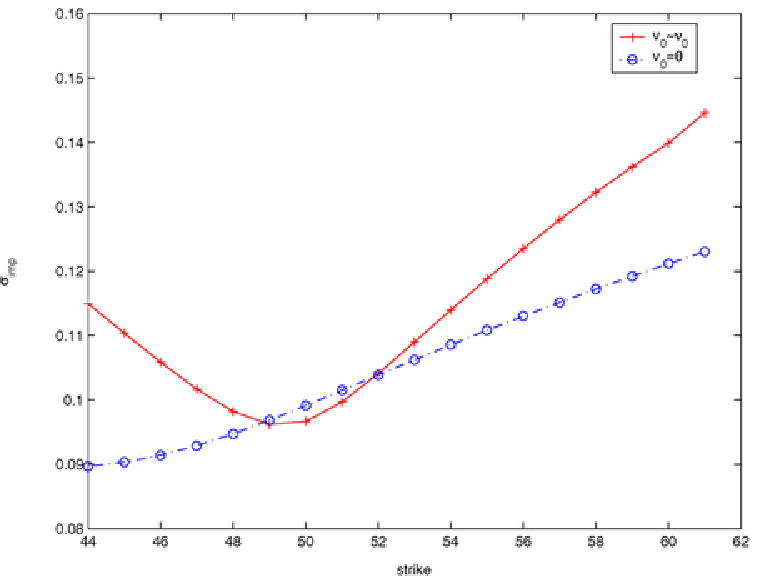}

\caption{$\rho=0.5$, $T\mapsto\sigma_{\mathrm{imp}}(50,T)$.}
\label{foncfigure6}
\vspace*{-14pt}
\end{figure}

The main differences between these models then appear for short
maturities. That is why we complete this part by a representation
of the volatility curve when $T=0.1$ for $\rho=0$ and $\rho=0.5$
in Figures \ref{foncfigure5} and \ref{foncfigure6}, respectively.
We observe that for short maturities, the volatility smile is more
curved and the skew is steeper. These phenomena seem
interesting for calibration since one well-known drawback of the
standard Heston model is that it can have overly
flat volatility curves for short maturities.

\begin{table}
\caption{Approximation of the Asian call price in the BNS model\label{callasian2}}
\begin{tabular*}{\textwidth}{@{\extracolsep{\fill}} lccccccc@{}}
\hline
$K$ & 44 &45& 46&47 & 48&49&50 \\
\hline
Asian call (ref.) & \textbf{6.75} &\textbf{5.83}& \textbf{4.93}& \textbf{4.05} & \textbf{3.18} & \textbf{2.35} &\textbf{1.57}\\
$N=5\cdot10^4$& 6.83 & 5.91& 5.01 &4.10& 3.22& 2.35& 1.51 \\
$N=5\cdot10^5$& 6.78 & 5.86& 4.96 & 4.06& 3.19 &2.34& 1.52 \\
$N=5\cdot10^4$ (CP parity)& 6.76 & 5.85& 4.94& 4.07& 3.20& 2.29& 1.51\\
$N=5\cdot10^5$ (CP parity)& 6.75 & 5.83& 4.93& 4.04& 3.17& 2.32& 1.54 
\end{tabular*}
\begin{tabular*}{\textwidth}{@{\extracolsep{\fill}} lcccccc@{}}
\hline
$K$ & 51 &52& 53&54 & 55&56 \\
\hline
Asian call (ref.)& \textbf{0.91}& \textbf{0.55}& \textbf{0.39} &\textbf{0.29} &\textbf{0.23} & \textbf{0.18} \\
$N=5\cdot10^4$& 0.77& 0.46& 0.33& 0.27& 0.22& 0.19\\
$N=5\cdot10^5$& 0.79 & 0.48& 0.34& 0.27& 0.21& 0.17 \\
$N=5\cdot10^4$ (CP parity)&0.79& 0.47& 0.37& 0.27& 0.23& 0.19\\
$N=5\cdot10^5$(CP parity)& 0.83& 0.50& 0.36& 0.28& 0.22& 0.17\\
\hline
\end{tabular*}
\vspace*{-6pt}
\end{table}

\subsection{Numerical tests on Asian options in the BNS SSV model}

The BNS model introduced in Barndorff-Nielsen and Shephard
\cite{barndorffshephard} is a
stochastic volatility model where the volatility process is a
L\'{e}vy-driven positive Ornstein--Uhlenbeck process. The dynamic of
the asset price $(S_t)$ is given by $S_t=S_0\exp(X_t)$,
\begin{eqnarray*}
\mathrm{d}X_t &=& \bigl(r-\tfrac{1}{2}v_t\bigr)\,\mathrm{d}t
+ \sqrt{v_t}\,\mathrm{d}W_t+\rho \,\mathrm{d}Z_t, \qquad\rho\le 0,
\\
\mathrm{d}v_t &=& -\mu v_t\,\mathrm{d}t+\mathrm{d}Z_t,\qquad\mu>0,
\end{eqnarray*}
where $(Z_t)$ is a subordinator without drift term and L\'{e}vy
measure $\pi$. In the following, we assume that $(Z_t)$ is a
tempered stable subordinator, that is, that
\[
\pi(\mathrm{d}y)=1_{\{y>0\}}
\frac{c\exp(-\lambda y)}{y^{1+\alpha}}\,\mathrm{d}y, \qquad
c>0, \lambda>0, \alpha\in(0,1).
\]
As in the Heston model, we want to use our algorithm as a way of
option pricing when the stochastic volatility evolves under its
stationary regime and test it on Asian options using the method
described in detail in Section \ref{methodpricing}. This model
does not require a specific discretization and the approximate
Euler scheme (P) (see Section \ref{levydriven}) relative to
$(v_t)$ can be implemented using the rejection method. In Table
\ref{callasian2}, we present our numerical results obtained for
the following choices of parameters, steps and weights:
\[
\rho=-1,\qquad
\lambda=\mu=1,\qquad
c=0.01,\qquad
\alpha=\tfrac{1}{2},\qquad
\gamma_n=\eta_n=n^{-{1}/{3}}.
\]
The computation times for $N=5.10^4$ and $N=5.10^5$ are about
8.5 s and 93 s, respectively.
Note that for this model, the convergence seems
to be slower because of the approximation of the jump component.

\section*{Acknowledgement}

The authors would like to thank Vlad Bally
for interesting comments on the paper.

\printhistory

\end{document}